\documentclass[final, 3p, compress,11pt,times]{elsarticle}
\usepackage{graphicx}
\usepackage{multirow}
\usepackage{framed}
\usepackage{amsmath}
\usepackage{bm}
\usepackage{amssymb}
\usepackage{framed}
\usepackage{psfrag}
\usepackage{epsfig}
\usepackage[usenames, dvipsnames]{xcolor}
\usepackage{array}
\usepackage{booktabs}
\usepackage[colorlinks=true]{hyperref}
\usepackage{tikz}
\usepackage{pgfplots}
\usepackage{diagbox}
\usepackage{enumitem}
\usepackage{lineno}
\usetikzlibrary{arrows}
\usepackage{caption}
\newcommand{\bzero}{\mathbf{0}}

\renewcommand{\div}{\mathrm{div}\,}

\newcommand{\bsg}{\bm{\sigma}}
\newcommand{\eps}{\varepsilon}
\newcommand{\bep}{\bm \eps}
\newcommand{\bI}{\mathbf{I}}

\newcommand{\dev}{\mathrm{dev}}

\newcommand{\tr}{\text{tr}\,}

\newcommand{\grad}{\text{grad}\,}

\newcommand{\be}{\begin{equation}}
\newcommand{\ee}{\end{equation}}
\newcommand{\ba}{\begin{array}}
\newcommand{\ea}{\end{array}}
\newcommand{\bve}{\boldsymbol{\varepsilon}}

\definecolor{cloudblue}{RGB}{162,190,227}
\DeclareMathAlphabet\mathbfcal{OMS}{cmsy}{b}{n}
\DeclareMathAlphabet\mathbfit{OML}{cmm}{b}{it}

\journal{}
\makeatletter
\def\ps@pprintTitle{%
  \let\@oddhead\@empty
  \let\@evenhead\@empty
  \let\@oddfoot\@empty
  \let\@evenfoot\@oddfoot
}
\makeatother

\begin{document}
\begin{frontmatter}
\title{Fast staggered schemes for the phase-field model of brittle fracture based on the fixed-stress concept}

\author[]{Chenyi Luo\corref{cor}}
\ead{cheluo@ethz.ch}
\cortext[cor]{Corresponding author.} 

\address[1]{Department of Mechanical and Process Engineering, ETH Z{\"u}rich, 8092 Z{\"u}rich, Switzerland}

\begin{abstract}
Phase field models are promising to tackle various fracture problems where a diffusive crack is introduced and modeled using the phase variable. Owing to the non-convexity of the energy functional, the derived partial differential equations are usually solved in a staggered manner. However, this method suffers from a low convergence rate, and a large number of staggered iterations are needed, especially at the fracture nucleation and propagation. In this study, we propose novel staggered schemes, which are inspired by the fixed-stress split scheme in poromechanics. By fixing the stress when solving the damage evolution, the displacement increment is expressed in terms of the increment of the phase variable. The relation between these two increments enables a prediction of the displacement and the active energy based on the increment of the phase variable. Thus, the maximum number of staggered iterations is reduced, and the computational efficiency is improved. We present three staggered schemes by fixing the first invariant, second invariant, or both invariants of the stress, denoted by S1, S2, and S3 schemes. The performance of the schemes is then verified by comparing with the standard staggered scheme through three benchmark examples, i.e., tensile, shear, and L-shape panel tests. The results exhibit that the force-displacement relations and the crack patterns computed using the fast schemes are consistent with the ones based on the standard staggered scheme. Moreover, the proposed S1 and S2 schemes can largely reduce the maximum number of staggered iterations and total CPU time in all benchmark tests. The S2 scheme performs comparably except in the shear test, where the underlying assumption is violated in the region close to the crack.

\end{abstract}
\end{frontmatter}

\section{Introduction}\label{intro}

The phase-field model has been applied to fracture problems for over twenty years. The pioneers, e.g., Marigo, Francfort, and Bourdin, founded its rationality from a mathematical point of view, such as $\Gamma$-convergence and regularization of the variational formulation, in their serial works \cite{francfort1998revisiting, bourdin2000numerical, bourdin2008variational}. The phase field model is competitive not only for its feasibility of simulating complex crack phenomena such as cracking nucleation, propagation, and branching but also for its simple numerical implementation, e.g., no change of geometry and no need to predict possible crack modes. For these advantages, the phase-field model has been immediately investigated by many groups from both the physics and mechanics communities. Till now, the benchmark achievements include splits of the strain energy for physically reasonable crack generation \cite{amor2009regularized,miehe2010phase,ambati2015review}, numerical implementation of the irreversibility \cite{bourdin2000numerical, miehe2010phase, gerasimov2019penalization}, coupling with other features, e.g., thermal cracking, hydraulic fracturing, plasticity \cite{bourdin2014morphogenesis,ehlers2017phase,miehe2016phase}, to name but a few.

Although the phase-field model has proven its potential in tackling complex fracture problems, its further development and application are impeded by computational costs. The costs mainly result from a fine mesh at least locally in the crack zone and the slow convergence in the numerical implementation. The first problem can be alleviated by applying the adaptive mesh refinement method \cite{heister2018parallel}. However, when a large number of cracks are distributed within the domain, e.g., thermal cracks or desiccation-induced cracks, a fine mesh will be required for the whole domain, which makes this method no longer efficient. To reduce the computational costs owing to the slow convergence, Gerasimov and De Lorenzis \cite{gerasimov2016line} proposed a line search assisted approach based on a monolithic scheme, where the coupling between the displacement and the crack is fully considered. Hence, the rate of convergence is improved. Nevertheless, the monolithic scheme usually leads to a larger Jacobi matrix and takes more time in one Newton-Raphson (NR) iteration. Another approach based on the staggered scheme and over-relaxation was developed \cite{storvik2021accelerated}, where the displacement increment is carefully exaggerated by a relaxation factor. The number of staggered iterations can be reduced if the relaxation factor is properly chosen. An analogous method was reported in \cite{brun2020iterative}, which introduced artificial friction terms to both the momentum balance and evolution equation. By varying the contribution of these two terms, fewer staggered iterations are required. However, for both methods, the choice of the numerical parameters dominates the computational efficiency and accuracy. As the latter method introduces errors to the system, the consistency between the results from the new schemes and the standard one cannot always be guaranteed.

In this context, we propose fast staggered schemes inspired by the fixed-stress split previously exploited in poromechanics \cite{kim2011stability,kim2011stability1,kim2011stability2}. The basic idea is to solve the evolution equation assuming that the stress is unchanged. Because if the evolution equation is solved, provided that the coupling term in the momentum balance, namely stress, keeps constant, the momentum balance will automatically hold since no terms are changed after solving the evolution equation. Note that if the phase variable changes when solving the evolution equation, the strain must also change, which results in a change of the active strain energy in the evolution equation. In our scheme, as the stress is unchanged, the update of the active strain energy will also be considered when solving the evolution equation. In this sense, our approach will achieve a higher convergence rate because we predict the active strain energy based on the change of the phase variable and the resultant strain change. In contrast, in the standard staggered scheme, this energy can only be updated after solving the momentum balance in the next staggered iteration. However, the challenge comes from deriving the relation between the strain and the phase variable to guarantee that the stress is constant. As the dimensions of the strain and phase variables are different, we propose to fix the first and second invariants of the stress instead of the stress itself. By this means, we are able to obtain the relation between the first and second invariants of strain and the phase variable. Thus, when updating the phase variable, a associated update is applied to the strain and the active energy such that the stress is unchanged and the momentum balance still holds.

The paper is organized as follows. Section \ref{sec:model} briefly introduces the phase-field model to brittle fracture providing the governing partial differential equations. In Section \ref{num-algo}, we present the standard staggered scheme and discuss its limitations. Then we set forth the derivation of the fast staggered schemes along with an alternative interpretation. Three numerical examples, including a tensile test, a shear test, and an L-shape panel test, are presented in Section \ref{sec:cases} to demonstrate the computational efficiency of the proposed numerical algorithms by comparing to the standard staggered scheme. A conclusive summary closes this study. Details on the derivation of the relations are included in \ref{summplementry_deduction}.

\section{Phase field model for brittle fracture}\label{sec:model}
In this section, the phase field model for brittle fracture is briefly introduced. The phase field model has been interpreted in different ways \cite{bourdin2000numerical, kuhn2010continuum,ehlers2017phase}. Here, we follow the idea of Bourdin, Francfort and Marigo \cite{bourdin2008variational}, among the earliest studies introducing the phase field model to fracture. The state of the material, either intact or cracked, is characterized by a phase variable, as
\be
d = 
\left\{
\ba{ll}
0, \quad & \text{for an intact solid},\\
1, \quad & \text{for a fully cracked solid}.
\ea
\right.
\ee
Thus, the energy functional for a pure solid subject to brittle fracture takes the form of 
\be
\mathcal{E} (\mathbf{u}, d) = \int_\Omega\underbrace{ \left[g(d) \Psi^+(\boldsymbol{\varepsilon} (\mathbf{u})) + \Psi^-(\boldsymbol{\varepsilon} (\mathbf{u}))\right]}_{\displaystyle \Psi} \text{d} v + \frac{G_c}{c_w} \int_{\Omega} \left[\frac{w(d)}{l} + l|\grad d|^2\right] \text{d}v- \int_{\partial \Omega} \mathbf{t} \cdot \mathbf{u} \text{d} s,
\ee
where the linearized strain tensor $\boldsymbol{\varepsilon}$ is defined by $\boldsymbol{\varepsilon} = \frac{1}{2}\left(\grad \mathbf{u}+\grad^T\mathbf{u}\right)$ with $\mathbf u=$ displacement.
The energy comprises strain energy, fracture energy, and work done by the external load $\mathbf{t}$. The strain energy is split into two parts to describe the preference for cracking under expansion or tensile stress. The positive or the active part, $\Psi^+$,  decreases to zero in a cracked material while the negative or the inactive part, $\Psi^-$, remains unchanged. Two prevalent splits are proposed by \cite{amor2009regularized} and \cite{miehe2010phase}, respectively. The work of \cite{ambati2015review} shows that both splits yield similar results in both the tensile and the shear tests. Although both splits have their difficulties in describing cracking under complex stress states \cite{luo2021phase, de2021nucleation},  we take the volumetric and deviatoric (VD) split from \cite{amor2009regularized} in the following discussion. The proposed staggered scheme in this paper should be easily extended to the phase field model with other splits, e.g., spectral decomposition.

When the VD split is applied, the active and inactive strain energies are given by
\be
\Psi^+ = \mu \, \bve^\text{dev} \cdot \bve^\text{dev} +\frac{k}{2} \langle\tr \bve\rangle_+^2, \quad \Psi^+ = \frac{k}{2} \langle\tr \bve\rangle_-^2,
\ee
where the subscripts are defined as
\be
   \left< \square \right>_+= \frac{ \square+|\square|}{2} \quad \text{and} \quad
   \left<\square\right>_-= \frac{\square -| \square|}{2},
\ee
and the Lam\'e constant and the bulk modulus are denoted by $\mu$ and $k$. Note that when $k$ appears as a superscript in the following, it denotes the $k$th staggered iteration.

\begin{table}[htbp]
\small
\caption{Parameters in AT1 and AT2 models}
\label{tab:at1at2}
	\centering\small
	\begin{tabular*}{0.5\textwidth}
			{@{\extracolsep\fill}cccc@{\extracolsep\fill}}
		\toprule
		\quad Model & $g(d)$& $w(d)$ & $c_w$\phantom{dd}\\
		\midrule
		\quad AT1 &\multirow{2}{*}{$(1-d)^2+\eta$}& $d$ & ${8}/{3}$\phantom{dd}\\\
		 \quad \hspace{-2pt}AT2 && \, $d^2$ & $ 2 $\phantom{dd} \\
		\bottomrule
	\end{tabular*}
\end{table}

The second part of $\Psi^+$ is the fracture energy which contains not only the phase variable, $d$, but also its gradient, $\grad d$. This gradient term makes the derived model no longer suffer from the pathological effect in the standard damage model, e.g., nonphysical dependence on the mesh size. Here, the critical energy release rate from Griffith's theory is denoted by $G_c$, and the crack width is denoted by $l$. The contribution of the strain energy and the fracture energy to the energy functional is weighted by the phase variable, or more precisely, by two functionals $g(d)$ and $w(d)$. Two common combinations of these two functions are known as the AT1 and AT2 models, as listed in Table \ref{tab:at1at2}, where AT comes from the initials of Ambrosio and Tortorelli \cite{ambrosio1990approximation}. The coefficient, $c_w$, is chosen to satisfy the $\Gamma-$convergence. It needs to emphasize the differences between the two AT models. As shown in \cite{pham2011gradient}, the phase variable in the AT1 model keeps zero before the strain energy exceeds a certain value. In contrast, the evolution of the phase variable in the AT2 model occurs even if the strain energy is very small. In other words, the AT1 model exhibits a purely elastic behavior before cracking, while the AT2 model has no such region. In this sense, the AT1 model is more suitable for describing brittle fractures in an elastic solid.

After minimizing the energy functional with respect to the primary variables, we derive the momentum balance equation in the absence of gravity and the evolution equation for the phase variable as
\be \label{eq_governing_equation}
\ba{c}
\div \boldsymbol{\sigma} = \mathbf{0},\\[5pt]
\displaystyle \frac{\partial g(d)}{\partial d} \Psi^+ + \frac{G_c}{c_w}\left[\frac{1}{l}\frac{\partial w(d)}{\partial d} - 2l \div \grad d\right]\geq 0.
\ea
\ee
The constitutive equation for the stress, $\bsg$, is given by
\begin{equation} \label{eq_sigma}
    \bsg=g(d)\left[
        2\mu\bep^{\dev}+k\langle\tr \bep\rangle_+\bI
    \right]
    +k\langle\tr \bep\rangle_-\bI.
\end{equation}
Note that the evolution equation is also subject to 
\be
\dot{d} \geq 0,
\ee
such that the Kuhn-Tucker condition is fulfilled here. Meanwhile, the temporal monotone of the phase variable, i.e., the irreversibility of the crack generation, is guaranteed.
It is worthwhile mentioning that the most popular and simple way to numerically realize the preceding inequality was proposed in \cite{miehe2010phase}, where the active strain energy is substituted by a history variable. The approach is widely adopted, particularly when the AT2 model is applied. However, this approach fails in the AT1 model, and more importantly, the model based on the history variable is no longer in a variational framework. In this context, we follow the idea from \cite{gerasimov2019penalization} and introduce a penalty term in the evolution equation,

\be \label{eq_evolution_equation}
\frac{\partial g(d)}{\partial d} \Psi^+ + \frac{G_c}{c_w}\left[\frac{1}{l}\frac{\partial w(d)}{\partial d} - 2l \div \grad d\right]  +  \gamma \left<d - d_{n-1}\right>_-=0,
\ee
where $\gamma$ is the penalty parameter and suggested to be 
\be
\gamma := \frac{27G_c}{64l \text{tol}_\text{ir}^2}.
\ee
Here, $\text{tol}_\text{ir}$ is user-prescribed irreversibility tolerance threshold and chosen as 0.01. Note that this method can apply to both the AT1 and AT2 models.

\section{Numerical algorithms}\label{num-algo}
In this section, we first summarize the standard staggered scheme and then propose novel staggered schemes, which are inspired by the fixed-stress split \cite{kim2011stability}. Comparative discussions on both schemes will be provided in Section \ref{sec_comment_standard} and \ref{sec_alternative_interpretation}, respectively.

\subsection{Standard staggered scheme}
\begin{figure}[!htp]
   \centering
   \includegraphics[width=\textwidth]{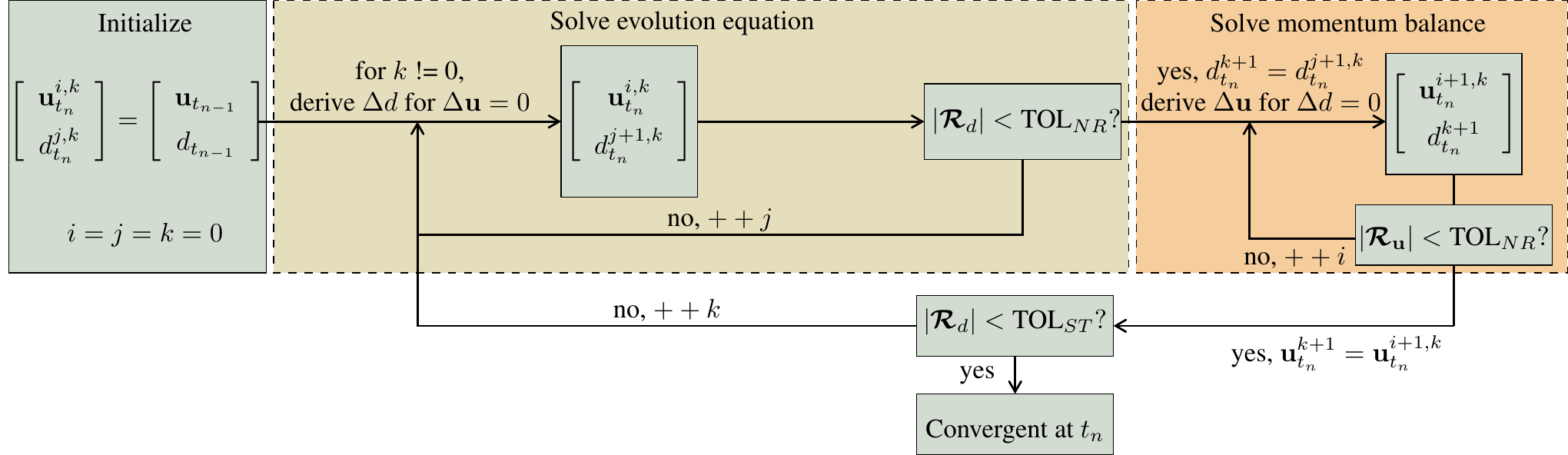}
   \caption{Flowchart of the standard staggered scheme.} 
   \label{fig:standard_staggered}
\end{figure}

The standard staggered scheme is a straightforward solution when one solves multiple equations in a sequential manner. The implementation is simple, i.e., solving one variable from its governing equation under the assumption that the other variables stay unchanged. For the phase-field model, two different staggered schemes are usually applied, and the difference lies in the sequence of solving equations. Here, an example based on the $d-\bf u$ sequence is illustrated in Fig. \ref{fig:standard_staggered}. At the time step $t_{n}$, the primary variables, $\mathbf{u}$ and $d$ are initialized with their values from the last time step $\mathbf{u}_{t_n}^{i,k}$ and $d_{t_n}^{j,k}$. The superscript ($i, \, j$) represents the iterator for the NR scheme solving $\mathbf{u}$ and $d$, respectively, while $k$ denotes the iteration step of the staggered scheme. After initiation, we first solve the evolution equation to get the new value for $d$. Two issues need to be mentioned here. First, this equation yields no update of $d$ for $k=0$ since the displacement update and the resultant energy change have not yet been considered. So, this step could be skipped for $k= 0$. Second, the evolution equation based on the AT1 model is nonlinear owing to the irreversibility of the phase variable and the penalty term, while the equation for the AT2 model is linear and can be solved with one iteration. After deriving $d_{t_n}^{k+1}$, we solve the momentum balance by fixing $d$ and yield $\mathbf{u}_{t_n}^{k+1}$. Then we check the residual of the evolution equation associated with $d_{t_n}^{k+1}$ and $\mathbf{u}_{t_n}^{k+1}$. If the residual exceeds the tolerance of the staggered scheme, we again solve the evolution equation and the momentum balance till the convergence is achieved.
It should be mentioned that two types of tolerance, namely the tolerance for the NR scheme and the one for the staggered scheme, usually follow the relation below
\be
\text{TOL}_{ST} > \text{TOL}_{NR}. 
\ee

\subsection{Discussions on the standard staggered scheme} \label{sec_comment_standard}
The readers may have already noticed the difference between a staggered scheme and a monolithic scheme for solving the phase field model. For an abstract discussion, we define the system as $\mathbf{F}$ which is composed of two parts, $\bf F_\mathbf{u}$ and $F_d$,
\be
\mathbf{F}=\left[\ba{c} \mathbf{F}_\mathbf{u}\\ \mathbf{F}_d\ea\right].
\ee 
These two functions are the weak forms of Eq. \eqref{eq_governing_equation}.1 and Eq. \eqref{eq_evolution_equation} after spatial discretization, which take the form as
\be
\mathbf{F}_\mathbf{u} = \int_\Omega \mathbf{B}_\mathbf{u}^T \bsg \text{d}v - \int_{\partial \Omega} \mathbf{N}_\mathbf{u}^T  \bar{\mathbf{t}} \cdot \text{d}\mathbf{a},
\ee
and 
\be
\mathbf{F}_d = \int_\Omega \left\{\mathbf{N}_d \left[-2 (1-d) \Psi^+ +  \frac{G_c}{c_w l}+ \gamma \left<d - d_{n-1}\right>_-\right]+2\mathbf{B}_d \frac{G_c l }{c_w} \grad d \, \right\}\text{d}v - \int_{\partial \Omega} \mathbf{N}_d \, \grad d \cdot\text{d} \mathbf{a} =\mathbf 0.
\ee
The shape functions for the displacement and the phase variable are represented by $\mathbf{N}_\mathbf{u}$ and $\mathbf{N}_d$, while the linearized strain operator and the derivative operator for the phase variables are denoted by $\mathbf{B}_{\mathbf{u}}$ and $\mathbf{B}_d$.
As the system is nonlinear, the increments for each NR iteration in a monolithic algorithm can be computed by solving the following linear system
\be \label{eq_monolithic_NR}
-\mathbf{R}=\left[\ba{c}-\mathbf{R}_\mathbf{u}\\-\mathbf{R}_d\ea\right] = \underbrace{
\left[
\ba{cc}
\mathbf{K}_{\mathbf{u} \mathbf{u}} &\mathbf{K}_{\mathbf{u} d} \\[5pt]
\mathbf{K}_{d \mathbf{u}} &{\bf K}_{d d} 
\ea
\right]}_{{\mathbf{K}}}
\left[\ba{c} \Delta \mathbf{u}\\ 
\Delta d\ea\right] ,
\ee
where $\mathbf{R}$ is the residual vector and $\mathbf{K}$ is the Jacobi matrix. The component of the Jacobi matrix, denoted by $\mathbf{K}_{\alpha \beta}$, is computed by
\be
\mathbf{K}_{\alpha \beta} = \frac{\partial \mathbf{F}_\alpha}{\partial \beta}, \quad \alpha, \,\beta \in \{\mathbf{u}, \, d\}.
\ee
If a staggered scheme with  the $d-\mathbf{u}$ sequence is applied, the off-diagonal term, $\mathbf{K}_{d \mathbf{u}}$, is neglected within one staggered iteration since the evolution equation is satisfied with fixed displacement. Then the increments are computed as
\be \label{eq_staggered_NR_st}
-\mathbf{R}=\left[\ba{c}-\bar{\mathbf{R}}_\mathbf{u}\\ -\bar{\bf R}_d\ea\right] = \underbrace{
\left[
\ba{cc}
\bar{\mathbf{K}}_{\mathbf{u} \mathbf{u}} &\bar{\mathbf{K}}_{\mathbf{u} d} \\[5pt]
\mathbf{0} &\bar{\bf K}_{d d} 
\ea
\right]}_{{\mathbf{K}}}
\left[\ba{c} \Delta \mathbf{u}\\ 
\Delta d\ea\right].
\ee
Here variables with an overbar, e.g., $\bar{\mathbf{R}}_\mathbf{u}$, indicate that they are evaluated by fixing the phase variable and an updated $\mathbf{u}$, while those without overbars (e.g., $\mathbf{R}_\mathbf{u}$) are obtained with both $d$ and $\bf u$ updated.
It should be mentioned that these two equations, \eqref{eq_monolithic_NR} and \eqref{eq_staggered_NR_st} are usually not equivalent. More importantly, it has been proven that such a staggered scheme for solving coupled problems does not always yield convergent results. Especially when the model is not set up in a variational framework, the Jacobi matrix is not symmetric. For example, in poromechanics, similar algorithms, known as the fixed-strain and drainage splitting schemes, are proven conditionally stable, cf.~\cite{kim2011stability1,kim2011stability2}. The standard staggered scheme applied to the phase-field model has proven its convergence in \cite{bourdin2007numerical}.  However, neglecting the coupling term leads to a slow convergence rate compared to monolithic schemes \cite{gerasimov2016line, lampron2021efficient}.

\subsection{Fast staggered schemes based on the tensor invariants}
To achieve a higher convergence rate without introducing any artificial term, we propose new staggered schemes by fixing the invariants of the coupling term. More specifically, if the coupling term in the momentum balance keeps constant when $d$ evolves in solving the evolution equation, the momentum balance simultaneously holds within the domain \footnote{Note the boundary condition for $\mathbf{u}$ has not been and cannot be considered till now.}. In the present model, the coupling term of the momentum balance is the stress $\bsg$. The stress increment between two staggered iterations is given by $\Delta \bsg$,
\begin{equation}
    \bsg^{k+1}=\bsg^{k}+\Delta\bsg,
\end{equation}
which, based on Eq. \eqref{eq_sigma}, can be computed by
\begin{equation} \label{eq_delta_stress}
\ba{l}
  \displaystyle  \Delta \bsg = \left\{
        g(d)[
            2\mu(\bI\otimes\bI)^{\overset{24}{T}}-\frac{2\mu}{3}(\bI\otimes\bI)
            +k\langle\bI\otimes\bI\rangle_+
        ]
        +k\langle\bI\otimes\bI\rangle_-
    \right\}\Delta \bep -  \\[12pt]
    \displaystyle  \phantom{ \Delta \bsg =}
    -2(1-d)\left[
        2\mu \bep^{\dev}+
        k<\tr \bep>_+\bI
    \right]\Delta d.
    \ea
\end{equation}
Here the variables with the subscript ``+" are defined as
\begin{equation} \label{eq_definition_+}
    I_+=\left\{
        \begin{array}{ll}
        1, &\tr \bep>0\\
        0, &\text{else}
        \end{array}
    \right., \quad
    \bI_+=\left\{
        \begin{array}{ll}
        \bI_+, &\tr \bep>0\\
        0, &\text{else}
        \end{array}
    \right.,\quad
    \langle\bI\otimes\bI\rangle_+=\left\{
        \begin{array}{ll}
        \bI\otimes\bI, &\tr \bep>0\\
        \bzero, &\text{else}
        \end{array}
    \right. .
\end{equation}
Now the challenge becomes deriving the relation between a second-order tensor $\Delta \bsg$ and a scalar $\Delta d$. We first compute the dot product between the stress increment and the identity tensor. Since $\Delta \bsg$ is zero, its trace is also zero,

\begin{equation} \label{eq_assumption1}
    \Delta \bsg \cdot \bI = \bzero.
\end{equation}
This leads to the relation between $\tr \Delta \bep$ and $\Delta d$ as
\begin{equation}\label{eq:trdep}
    \tr \Delta\bep = \frac{2(1-d)}{g(d)}\langle\tr \bep\rangle_+\Delta d.
\end{equation}
Then we take the dot product between $\Delta \bsg$ and $\bep$, or between $\Delta \bsg$ and $\bep^\text{dev}$,
\begin{equation} \label{eq_assumption2}
    \Delta \bsg \cdot \bep^{\dev}
    =0 \quad (\Delta\bsg\cdot\bep = 0),
\end{equation}
and yield the second relation
\begin{equation}\label{eq:epdep}
    \bep^{\dev}\cdot \Delta \bep
    =\frac{2(1-d)}{g(d)}
    (\bep^{\dev} \cdot \bep^{\dev})
    \Delta d.
\end{equation}
Then we recall the evolution equation \eqref{eq_monolithic_NR}
\begin{equation} \label{eq_rd_relation}
    -\mathbf{R}_d =\mathbf{K}_{d \mathbf{u}} \cdot \Delta \mathbf{u}+{\bf K}_{dd}\Delta d = \mathbf{K}_{d \bep} \cdot \Delta \bep + {\bf K}_{dd} \Delta d
    \end{equation}
    and find out that the first term can be rewritten by
\begin{equation}
\mathbf{K}_{d \bep} \cdot \Delta \bep  =- 2 \mathbf{N}_d  (1-d) \, \boldsymbol{\sigma} \cdot \Delta \bep = -2\mathbf{N}_d (1-d) (2\mu\bep^\dev \Delta\bep+k\langle\tr \bep\rangle_+\tr\Delta\bep) = \tilde{\bf K}_{d\bep}  \Delta d,
\end{equation}
where
\begin{equation}\label{eq:dpsi}
   \tilde{\bf K}_{d\bep} = - 4\mathbf{N}_d \frac{(1-d)^2}{g(d)} [2\mu 
    (\bep^{\dev} \cdot \bep^{\dev})
+k \langle\tr \bep\rangle_+^2].
\end{equation}
Inserting Eq. \eqref{eq:dpsi} into Eq. \eqref{eq_rd_relation}, we derive the increments of the phase variable provided that parts of the stress, $\Delta \bsg \cdot \bI$ and $\Delta \bsg \cdot \bep$,  keep unchanged. 
In the next step, we need to update the coupling term, i.e., $\Psi^+$, in the evolution equation\footnote{The update of the coupling term only occurs after the first NR iteration of the evolution equation. For further iterations, this energy keeps unchanged. The reasons behind this are discussed as follows.
\begin{itemize}
\item As the boundary condition of $\mathbf{u}$ is not considered here, the evolution equation under zero stress increments is solved such that the reaction force on the boundary is unchanged. In other words, only the Neumann boundary condition is applied to the displacement when solving the evolution equation, which may result in divergent results. Especially for a fracture problem, continuous updates of the active strain energy lead to unlimited crack propagation till the crack goes through the whole domain. Meanwhile, the material is still under constant reaction force.
\item As already known, the energy functional with respect to ($\mathbf{u}$, $d$) is non-convex while it is convex when fixing either $\mathbf{u}$ or $d$. If the coupling term is considered at each NR iteration, the convexity might be violated, leading to divergent results.
\end{itemize}}. This active strain energy can be computed by 
\begin{equation}
    \Psi^{+,k+\frac{1}{2}}=
    \mu \, \bep^{\dev,k+\frac{1}{2}}\cdot\bep^{\dev,k+\frac{1}{2}}+\frac{k}{2}(\langle\tr \bep\rangle^{k+\frac{1}{2}}_+)^2,
\end{equation}
where the superscript $\frac{1}{2}$ represents the half step between staggered iteration $k$ and $k+1$, which is computed by the assumptions given in Eq. \eqref{eq_assumption1} and Eq. \eqref{eq_assumption2} but has not been adapted to the boundary condition of $\mathbf{u}$.
With Eq. \eqref{eq:trdep} at hand, the new strain trace can be computed via
\begin{equation}\label{eq:epn1}
    \langle\tr\bep\rangle_+^{k+\frac{1}{2}}=
    \langle\tr\bep\rangle_+^{k} +
    \langle\Delta \tr\bep\rangle_{+} =  \langle\tr\bep^k\rangle_+ +  \frac{2(1-d^k)}{g(d^k)}\langle\tr \bep^k\rangle_+\Delta d.
\end{equation}
For the sake of conciseness, the derivation of $\bep^{\dev,k+\frac{1}{2}}\cdot\bep^{\dev,k+\frac{1}{2}}$ is not shown here but provided in \ref{summplementry_deduction}. The deduction leads to
\be \label{e_dev_1/2}
\bep^{\dev,k+\frac{1}{2}}\cdot\bep^{\dev,k+\frac{1}{2}} = \bep^{\dev,k}\cdot \bep^{\dev,k} + 4    (\bep^{\dev,k} \cdot \bep^{\dev,k}) \left[\frac{(1-d^k)}{g(d^k)} \Delta d + \frac{(1-d^k)^2}{g^2(d^k)} (\Delta d)^2\right].
\ee
With Eq. \eqref{eq:epn1} and Eq. \eqref{e_dev_1/2}, we derive the new active strain energy, based on which the evolution equation is solved in a staggered manner. If the AT2 model, together with the history variable, is adopted, the energy update is no longer needed since the equation is linear and solved within one timestep. 

\be
\mathbf{u}_{t_n}^{i+1,k} = \mathbf{u}_{t_n}^{i,k} + \Delta \mathbf{u}.
\ee

\subsection{Alternative interpretation of the proposed staggered scheme} \label{sec_alternative_interpretation}
In the above section, we have derived the relation between the increments of $\mathbf{u}$ and $d$ from the momentum balance by taking the dot product of the stress increment with either the identity tensor or the strain tensor. As a result, instead of the relation between a second-order tensor and a scalar, we derive the relation between several scalar-valued quantities of a tensor and a scalar. Here, we attempt to interpret the above derived relations, i.e., Eq. \eqref{eq:trdep} and Eq. \eqref{eq:epdep} from another point of view. Recall that principal invariants of a tensor keep constant if the tensor is unchanged. Increments of the first and second principal invariants of the stress tensor yield
\be \label{eq_no_idea_for_name}
\frac{\text{d} }{\text{d} t} I_1^{\bsg} = \Delta \bsg \cdot \bI \quad \text{ and } \quad \frac{\text{d} }{\text{d} t} I_2^{\bsg} =- \bsg^\text{dev} \cdot \Delta \bsg.
\ee
Expanding the second equation, we can write the increment of the second principal invariants into two parts as
\be \label{eq_no_idea_for_name2}
\frac{\text{d} }{\text{d} t} I_2^{\bsg} =-[g(d) k \langle\tr \bep\rangle_+ + k\langle\tr \bep\rangle_-]\Delta \bsg \cdot \bI - g(d) (2 \mu \Delta \bsg \cdot \bep^\text{dev}).
\ee
The first term of the right-hand side of the preceding equation contains $\Delta \bsg \cdot \mathbf{I}$ while the second has $\Delta \bsg \cdot \bep^\text{dev}$. Comparing Eq. \eqref{eq_no_idea_for_name}, Eq. \eqref{eq_no_idea_for_name2} with Eq. \eqref{eq:trdep} and Eq. \eqref{eq:epdep}, we propose three staggered schemes, denoted by S1, S2, and S3, respectively, where three assumptions are respectively applied in solving the evolution equation. In the S1 scheme, we only consider the relation in Eq. \eqref{eq:trdep}, which is equivalent to an assumption of a constant first principal invariant of the stress tensor given by Eq. \eqref{eq_no_idea_for_name}. In the S2 scheme, the relation in Eq. \eqref{eq:epdep} is applied. Compare Eq. \eqref{eq_no_idea_for_name2}, this relation holds for a constant second invariant of the stress tensor with a negative strain trace. The last scheme, S3, which involves both Eq. \eqref{eq:trdep} and Eq. \eqref{eq:epdep}, guarantees that both the first and second invariants are fixed during the evolution of the phase variable. Table \ref{tab_summary_staggered_scheme} compares the standard and the proposed staggered schemes.
\begin{table}[htbp]
\small
\caption{Summary of the assumptions applied in the standard and the proposed staggered schemes.}
\label{tab_summary_staggered_scheme}
 \centering\small
 \begin{tabular*}{\textwidth}
   {@{\extracolsep\fill}ccccc@{\extracolsep\fill}}
  \toprule
  \vspace{1pt}
  Scheme  & ST & S1& S2 &S3 \\
  \midrule
  &&&&\\
  \multirow{3}{*}{Assumptions} &\multirow{3}{*}{$ \mathbf{u} = 0$}  & $ \Delta \bsg \cdot \bI   = 0 $ & $\Delta \bsg \cdot \bep^{\dev}=0 $ & $\Delta \bsg \cdot \bI   = 0 $  and $\Delta \bsg \cdot \bep^{\dev} = 0$ \\[8pt]
   & & $ \displaystyle \frac{\partial }{\partial t} I_1^{\bsg} =0$ & $\displaystyle \frac{\partial }{\partial t} I_2^{\bsg} =0$ and $\tr \bep < 0$&$\displaystyle  \frac{\partial }{\partial t} I_1^{\bsg} =0$  and $ \displaystyle \frac{\partial }{\partial t} I_2^{\bsg} =0$  \\[2pt]
   &&&&\\
  \bottomrule 
 \end{tabular*}
\end{table}
Recalling Eq. \eqref{eq_monolithic_NR} and Eq. \eqref{eq_staggered_NR}, the increments solved by the proposed staggered scheme fulfill
\be \label{eq_staggered_NR}
-\mathbf{R}=\left[\ba{c}-\tilde{\mathbf{R}}_\mathbf{u}\\ -\tilde{R}_d\ea\right] = \underbrace{
\left[
\ba{cc}
\tilde{\mathbf{K}}_{\mathbf{u} \mathbf{u}} &\tilde{\mathbf{K}}_{\mathbf{u} d} \\[5pt]
\mathbf{0} &{\bf K}_{d d} + \tilde{\bf K}^S_{d\bep}
\ea
\right]}_{{\mathbf{K}}}
\left[\ba{c} \Delta \mathbf{u}\\ 
\Delta d\ea\right].
\ee
The extra stiffness component $\tilde{\bf K}^S_{d\bep}$ may take three forms according to Table \ref{tab_summary_staggered_scheme}, i.e., 
\be
\ba{l}
\displaystyle \tilde{\bf K}^{S1}_{d\bep}  =   - 4\mathbf{N}_d \frac{(1-d)^2}{g(d)}
k \langle\tr \bep\rangle_+^2, \\[8pt]
\displaystyle \tilde{\bf K}^{S2}_{d\bep}  =  - 4\mathbf{N}_d \frac{(1-d)^2}{g(d)}(2\mu 
    \, \bep^{\dev} \cdot \bep^{\dev}),\\[8pt]
\displaystyle \tilde{\bf K}^{S3}_{d\bep}  =   - 4\mathbf{N}_d \frac{(1-d)^2}{g(d)} [2\mu 
    (\bep^{\dev} \cdot \bep^{\dev})
+k \langle\tr \bep\rangle_+^2].
\ea
\ee

\subsection{Numerical implementation}
 The numerical implementation of the proposed staggered scheme is rather convenient. Based on the standard staggered scheme, the extra stiffness components $\tilde{\bf K}_{d\bep} $ need to be computed at the first NR iteration. Then based on the derived $\Delta d$, the active strain energy $\Psi ^{+, k+\frac{1}{2}}$ is then updated by Eq. \eqref{eq:epn1} and Eq. \eqref{e_dev_1/2}. Note that both steps only occur when solving the evolution equation.
 
 \begin{figure}[!htp]
   \centering
   \includegraphics[width=\textwidth]{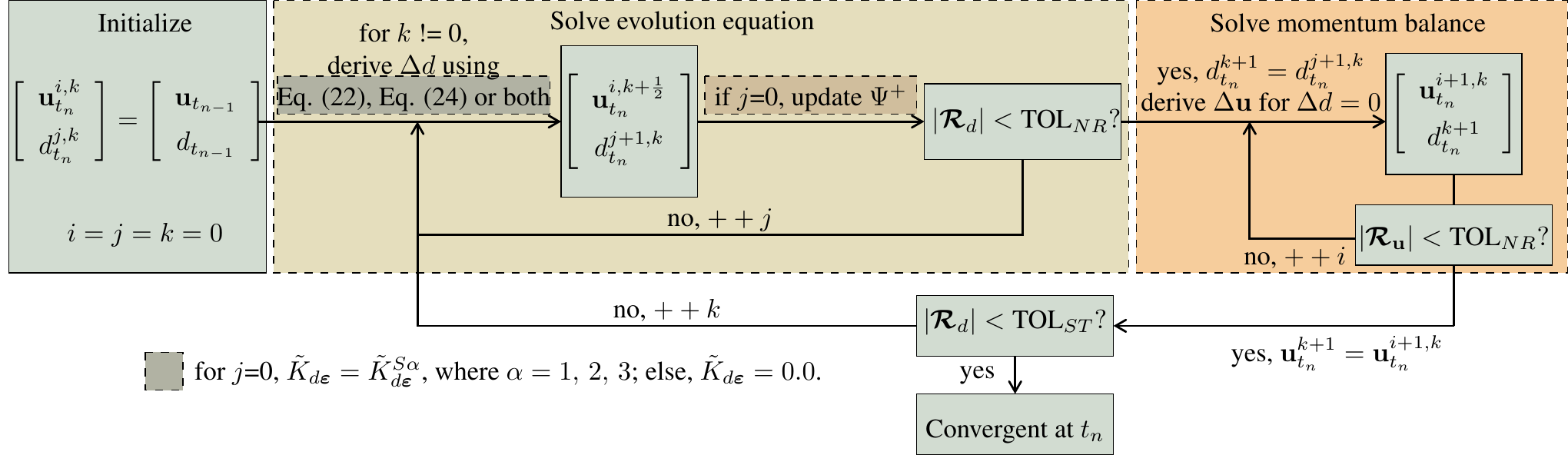}
   \caption{Flowchart of the staggered schemes based on the fixed-stress concept.} 
   \label{fig:standard_staggered_optimized}
\end{figure}
\begin{figure}[!htp]
   \centering
   \includegraphics[width=\textwidth]{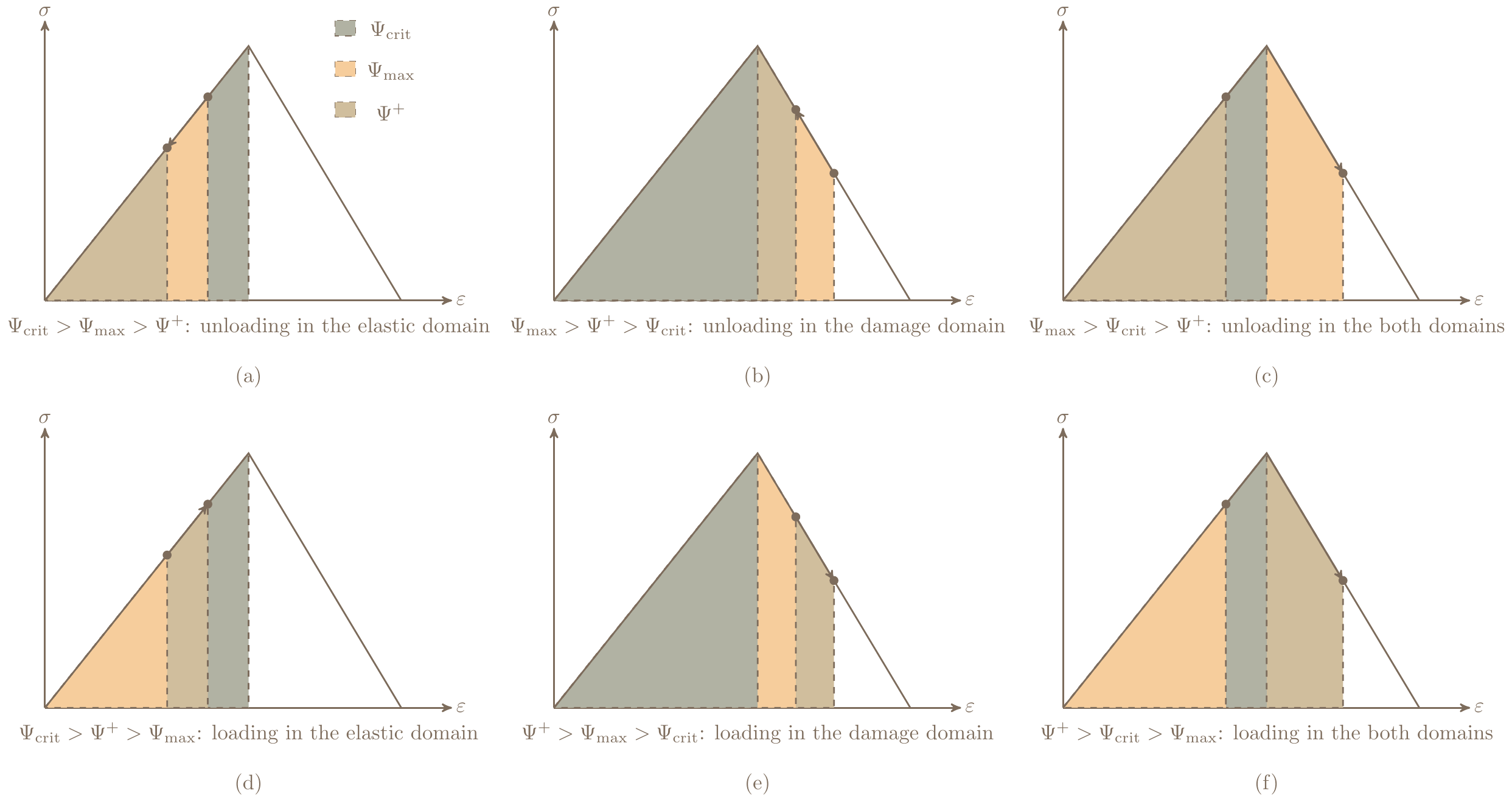}
   \caption{Illustration of the energy functionals, $\Psi_{\text{crit}}$, $\Psi^+$ and $\Psi_{\text{max}}$ under different loading and unloading conditions.} 
   \label{fig_loading_unloading}
\end{figure}

 However, attention needs to be paid to the fact that the phase variable severely deviates from its solution after the first NR iteration for the domain $d < d_{t_{n-1}}$. This phenomenon is because the penalty term has not been activated. Moreover, for $\dot{d} =0$, which corresponds to either an elastic behaviour or an unloading case, the stress is constant only for $\dot{\bep}=\mathbf{0} $. Therefore, the extra stiffness component and the update of the active strain energy will only be considered in the region where softening occurs. In the numerical implementation, it is realized by
 \be\label{condition_update}
 \tilde{\bf K}_{d\bep} =  \left\{
 \ba{ll}
  \tilde{\bf K}^{S \alpha }_{d\bep} \text{ with } \alpha =1 ,\, 2, \, 3 \quad& \text{for }\Psi^+ > \Psi_\text{crit} , \Psi^+ >\Psi_\text{max} \text{ and } \phi < 0.95\\
  0 & \text{else},
 \ea
 \right.
 \ee
where $\Psi_\text{crit} = \frac{G_c}{c_w l}$ and $\Psi_\text{max} = \text{Max} \, \Psi^+$. The underlying criteria can be explained by Fig. \ref{fig_loading_unloading}. When the current active energy $\Psi^+$ is smaller than $\Psi_\text{max}$, the material is unloading, cf. Fig. \ref{fig_loading_unloading}(a-c), where the phase variable does not evolve. When this energy is smaller than $\Psi_\text{crit}$, cf. Fig. \ref{fig_loading_unloading}(a) and (d), the material is in an elastic domain, where phase variable will neither evolve. Only when $\Psi^+ > \Psi_\text{crit} $ and $ \Psi^+ >\Psi_\text{max}$, the material is in the softening region and the phase variable is increasing, which corresponds to Fig. \ref{fig_loading_unloading} (e-f). Moreover, when the phase variable reaches its theoretical maximum 1, both the stiffness component $\tilde{\bf K}_{d\bep}$ and the increment of the damage variable $\Delta d$ are reduced to zero, and $\Psi^+$ will not be updated anymore. However, the real value of $d$ in numerical computation always varies around 1, which causes unnecessary updates of the active strain energy. Therefore, we set a limit of 0.95. When the phase variable is above 0.95, we consider that the material is fully broken, and no further update of the strain energy is needed.

After determining the stiffness components and deriving $\Delta d$, one needs to update the active strain energy. Note $\dot{d}\geq 0$ guarantees that the active strain energy increases monotonically during the update. However, $\Delta d$ might be a negative value in the numerical implementation. Therefore, when updating the active strain energy, Eq.~\eqref{eq:trdep} and Eq.~\eqref{eq:epdep} are modified as
\be
    \tr \Delta\bep = \left\{
    \ba{ll}
  \displaystyle  \frac{2(1-d)}{g(d)}\langle\tr \bep\rangle_+\Delta d, \quad& \text{for }\Delta d > 0,\\[8pt]
    0,&\text{else} ,
    \ea
    \right.
\ee
and
\begin{equation}
    \bep^{\dev}\cdot \Delta \bep
    =\left\{
    \ba{ll}
  \displaystyle  \frac{2(1-d)}{g(d)}
    (\bep^{\dev} \cdot \bep^{\dev})
    \Delta d,\quad &\text{for }\Delta d > 0,\\[8pt]
    0,&\text{else} .
    \ea
    \right.
\end{equation}

\section{Numerical examples}\label{sec:cases}
Three numerical examples are set up to verify the performance of the proposed staggered schemes. The tolerances for all examples are set as
\be
\text{TOL}_{NR} = 10^{-7} \quad \text{and}\quad \text{TOL}_{ST} = 10^{-6}.
\ee
Three-dimensional simulations are performed in all the cases, with the sample thickness chosen as the width of one element. Moreover, the movement in the thickness direction is forbidden such that the plain strain condition is fulfilled in the three-dimensional setting. As the tasks are computational-intensive, parallel computation is exploited, and the computation code is adapted to the framework of deal.ii \cite{arndt2021deal}. The computation was performed on 128 cores (for tensile and shear tests) and 256 cores (for the L-shape panel test) of the Euler cluster at ETH Z\"urich.

\subsection{Tensile tests}\label{sec:tensile}

\begin{figure}[!htp] 
   \centering
   \includegraphics[width=.9\textwidth]{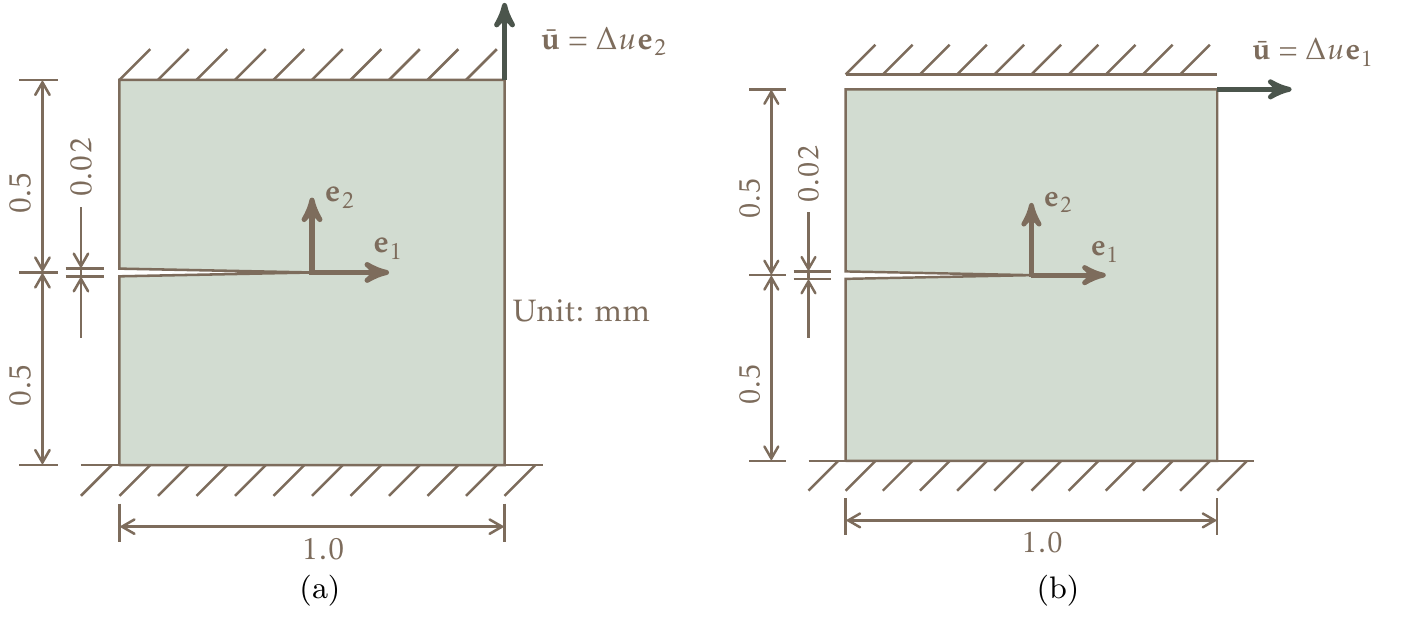}
   \caption{The geometry and the boundary conditions of the benchmark examples: (a) the tensile test; (b) the shear test.} 
   \label{fig_geo_bc_tensile}
\end{figure}

The first example is a benchmark model, which has been widely studied in the literature \cite{miehe2010phase, ambati2015review, wu2017unified}. We take the parameters originated from \cite{miehe2010phase}, as listed in Table \ref{tab:Parameters for the tensile and shear tests}. The same parameter setting will be  applied in the subsequent shear test. The geometry and the boundary conditions are illustrated in Fig. \ref{fig_geo_bc_tensile}, where the sample bottom is fixed, and the top boundary is subject to a linearly increasing vertical displacement. The loading step is given as
\be
\Delta u = 3\times 10^{-5} \text{ mm/s}, \qquad \Delta t_n=\left\{
\ba{ll}
170, \quad&t \leq 170 \text{ s},\\
1,\quad& t \leq 220 \text{ s}.
\ea
\right.
\ee

\begin{table}[htp]
\centering
\small
\vspace{5mm}
\caption{Parameter setting for the tensile and shear tests}
\label{tab:Parameters for the tensile and shear tests}       
\begin{tabular}{cccccc}
  \toprule
&&&&&\\[-4pt]
Parameter & $\mu$ & $l$ & $\lambda$  &  $\eta$ &$G_c$  \\[4pt]
  \midrule
&&&&&\\[-4pt]
Value & $8.077 \times 10^{4} \, \text{MPa}$ & $0.01\,\text{mm}$& $2.019 \times 10^{5}\, \text{MPa}$ &  $1\times 10^{-6}$ & $2.7 \,\text{N/mm}$ \\[4pt]
  \bottomrule
\end{tabular}
\end{table}

\begin{figure}[!htp]
   \centering
   \includegraphics[width=\textwidth]{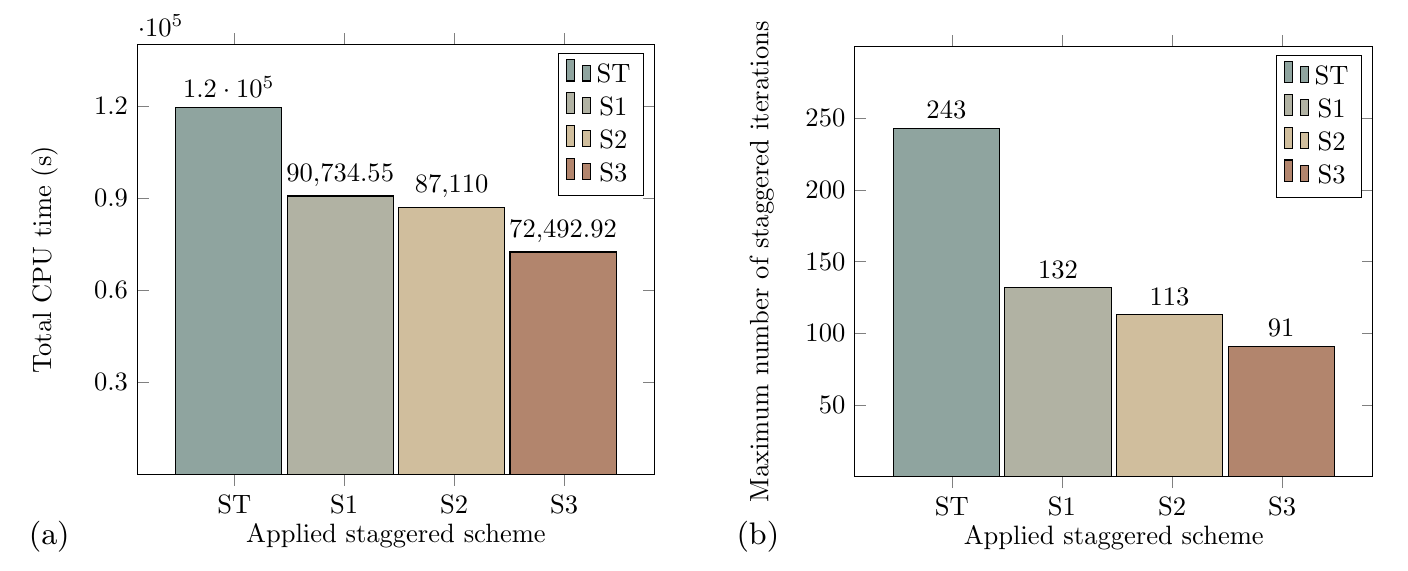} 
   \caption{Computational costs in the tensile tests: (a) total CPU time; (b) the maximum number of iterations for one timestep.}
   \label{fig_tensile_summary}
\end{figure}
When using the staggered schemes listed in Table \ref{tab_summary_staggered_scheme}, consistent results and cracking patterns are obtained, including
\begin{itemize}
\item A linear relation between the displacement and the vertical reaction force is found before cracking. After cracking, a sudden drop in the force is observed
\item A horizontal crack through the sample at the middle height is formed within one timestep.
\item The maximum number of staggered iterations occurs when the crack is initialized.
\end{itemize}
\begin{figure}[!htp]
   \centering
   \includegraphics[width=\textwidth]{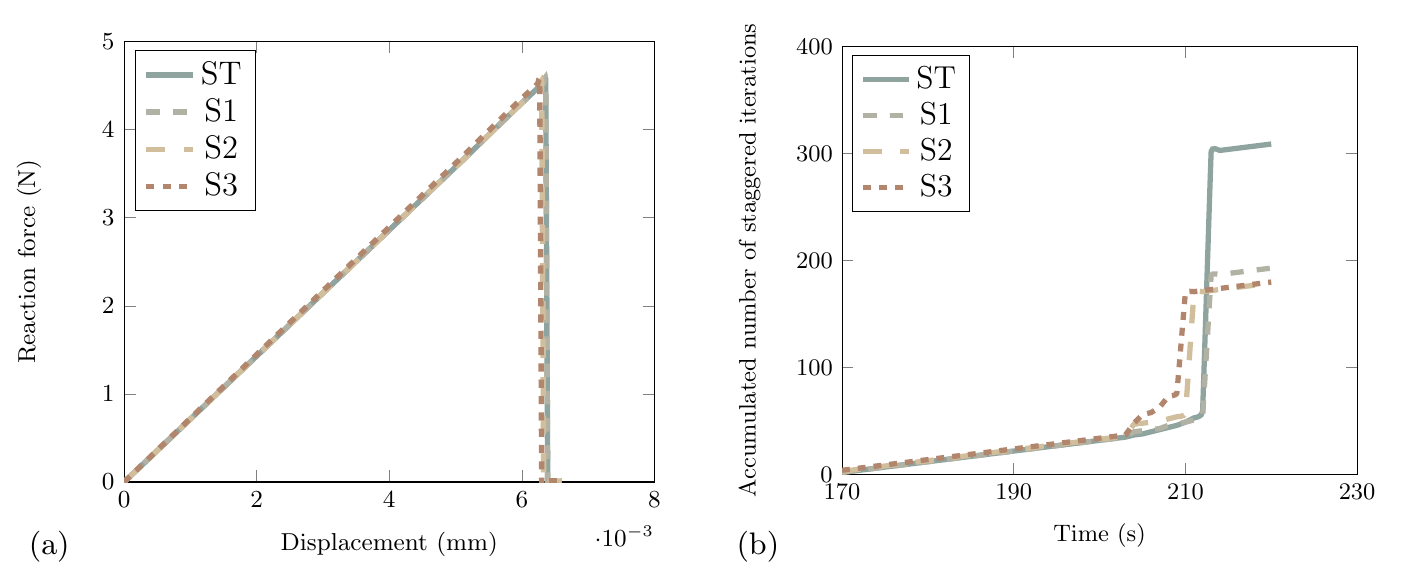}
   \caption{Results in the tensile tests: (a) force-displacement curves; (b) accumulated number of staggered iterations.}
   \label{fig_displacement_accu_iter_tensile}
\end{figure}
\begin{figure}[!htp]    \centering
   \includegraphics[width=\textwidth]{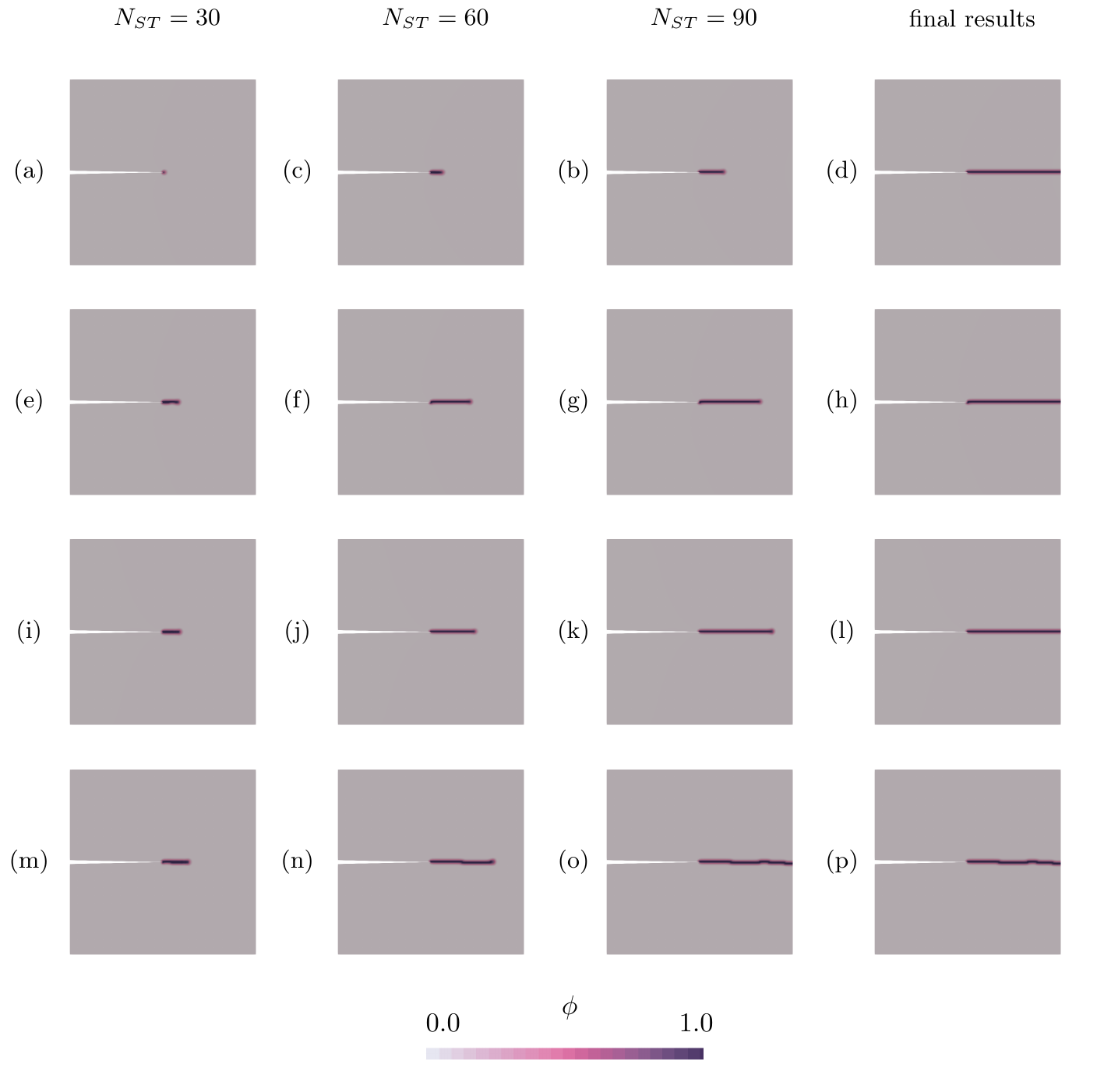}
   \caption{Crack propagation within the staggered iterations in the tensile test obtained using different staggered schemes: (a-d) ST scheme at 213 s; (e-h) S1 scheme at 213 s; (i-l) S2 scheme at 211 s; (m-p) S3 scheme at 210 s.}
   \label{fig_crack_propagation_tensile}
\end{figure}

Fig. \ref{fig_displacement_accu_iter_tensile}(a) compares the relations between the reaction force and the forced displacement. Fig. \ref{fig_crack_propagation_tensile} illustrates the crack evolution within one time step at different staggered iteration steps. In general, the results of force-displacement relations and the crack evolution obtained from different schemes agree well. The only difference is seen in Fig. \ref{fig_crack_propagation_tensile}(p), where the through crack bends slightly downwards at the final stage when using the S3 scheme. By a more careful comparison of the results in Fig. \ref{fig_tensile_summary}, Fig. \ref{fig_crack_propagation_tensile} and Fig. \ref{fig_displacement_accu_iter_tensile}(b), we perceive the following differences.
\begin{itemize}
\item[$\blacktriangleright$] When the S2 and S3 schemes are applied, cracking occurs slightly earlier. For the ST and S1 schemes, cracking occurs at time 213 s, while for the S2 and S3 schemes, it happens at 211 s and 210 s, respectively.
\item[$\blacktriangleright$] The total CPU time varies among the cases using different staggered schemes and is ranked as $T_\text{cpu}^{ST} > T_\text{cpu}^{S1}>T_\text{cpu}^{S2}>T_\text{cpu}^{S3}$, with the superscript indicating the adopted scheme. Using the S3 scheme reduces the CPU time to 60\% of that consumed in using the ST scheme. It is noted that CPU time cannot directly reflect the computation efficiency. \footnote{Although the CPU time measures the time consumed during the computation, it is not always reliable to measure the computational efficiency solely by this value for the below reasons:
\begin{itemize}
\item The CPU time is not a fixed value, which is randomly distributed among a certain range for computation on a single core.
\item When the parallel computation is applied, the variation owing to memory states and communication introduces more randomness to the actual CPU time.
\item When the parallel computation is conducted on a cluster, more variables, e. g., core distribution and different CPU types, increase the variation of the CPU time.
\end{itemize}
For this paper, we limit not only the total number of cores but also the number of cores per node and the CPU type, which helps to reduce the variation of the CPU time. Thereby, the variation of the CPU time is within 5\%. Compared to the CPU time, the staggered iterations are a stable value that only slightly fluctuates within 1. Hence, it is recommend to take both indicators together in evaluation of the computational efficiency.
}
\item[$\blacktriangleright$] The maximum number of iterations corresponding to different staggered schemes is proportional to the total CPU time, ranked as $N_\text{ST, max}^{ST} > N_\text{ST, max}^{S1}>N_\text{ST, max}^{S2}>N_\text{ST, max}^{S3}$. The maximum numbers of iterations corresponding to the S1, S2 and S3 staggered scheme are respectively 54\%, 46.5\% and 37\% of $N_\text{ST, max}^{ST}$ corresponding to the ST scheme.
\item[$\blacktriangleright$] The difference of the final accumulated staggered iterations among four staggered schemes is traced back to the difference in the maximum iteration number of one timestep. Before and after cracking, convergent results are usually achieved within one staggered iteration.
\end{itemize}

By comparing the crack patterns at the same staggered iteration step, e.g., Fig. \ref{fig_crack_propagation_tensile}(a), \ref{fig_crack_propagation_tensile}(e), \ref{fig_crack_propagation_tensile}(i) and \ref{fig_crack_propagation_tensile}(m), it can be seen  that by updating the active strain energy based on coupling relation, the crack propagates faster as compared to the ST model where this coupling term is completely ignored. As the S3 scheme fully considers both parts of the active strain energy, i.e., the parts induced by the deviatoric strain and the positive strain trace, the crack propagates most rapidly.

\subsection{Shear tests}\label{sec:shear}

The second example is also a benchmark model frequently used to test the performance of different energy split algorithms for phase-field fracture models, e.g., \cite{amor2009regularized, miehe2010phase, ambati2015review}. The material parameter setting is identical to that in the first example (Table \ref{tab:Parameters for the tensile and shear tests}). Fig. \ref{fig_geo_bc_tensile}(b) illustrates the geometry and boundary conditions. Unlike the tensile test, the sample is subject to a horizontal displacement on the top boundary. The loading step is given by
\be
\Delta u = 3\times 10^{-5} \text{ mm/s}, \qquad \Delta t_n=\left\{
\ba{ll}
250, \quad&t \leq 250 \text{ s},\\
1,\quad& t \leq 380 \text{ s}.
\ea
\right.
\ee

\begin{figure}[!htp]
   \centering
   \includegraphics[width=\textwidth]{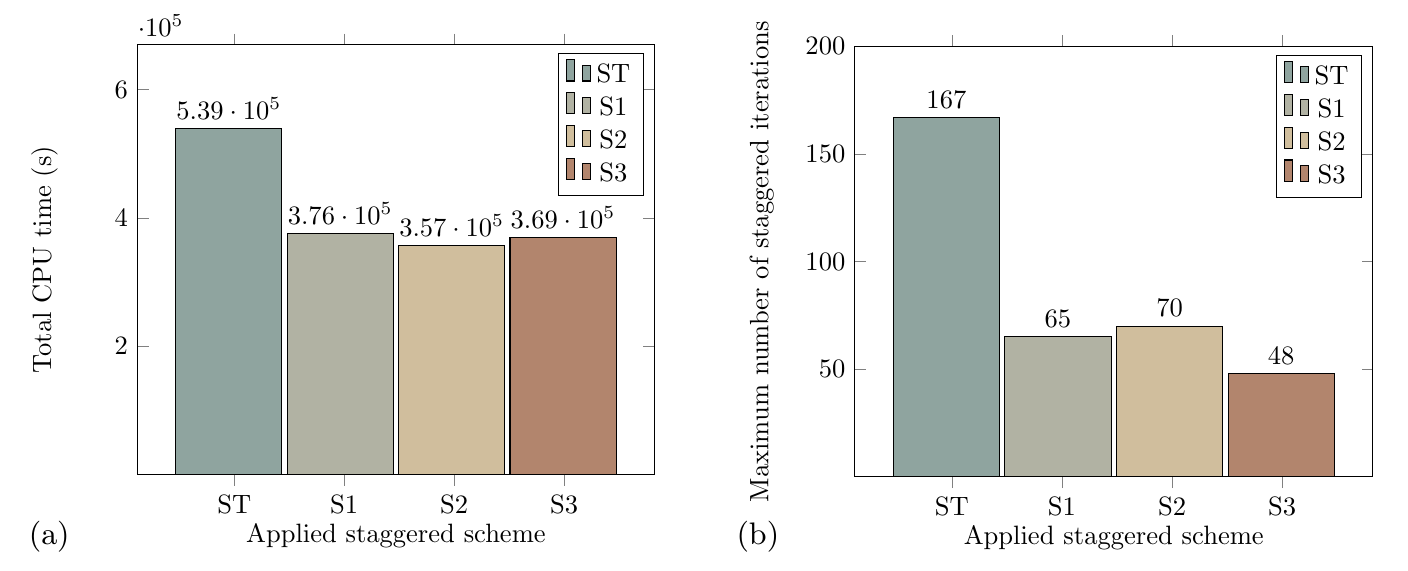}
   \caption{Comparison of the computational cost in the shear test: (a) total CPU time; (b) the maximum number of iterations for one timestep.}
    \label{fig_shear_summary}
\end{figure}

Comparing the results derived by different staggered schemes, we find they are pretty consistent, as discussed below.
\begin{itemize}
\item A linear relation between the displacement and the vertical reaction force is fulfilled before crack nucleation.
\item During the formation of a crack tip, the force has slightly oscillated under an increasing displacement.
\item After the crack length exceeds a certain value, the force decreases with an increasing displacement.
\item The crack is not formed within one timestep and the maximum number of iterations does not always occur at the crack nucleation but can take place in the softening region.
\end{itemize}
\begin{figure}[!htp] 
   \centering
   \includegraphics[width=0.8\textwidth]{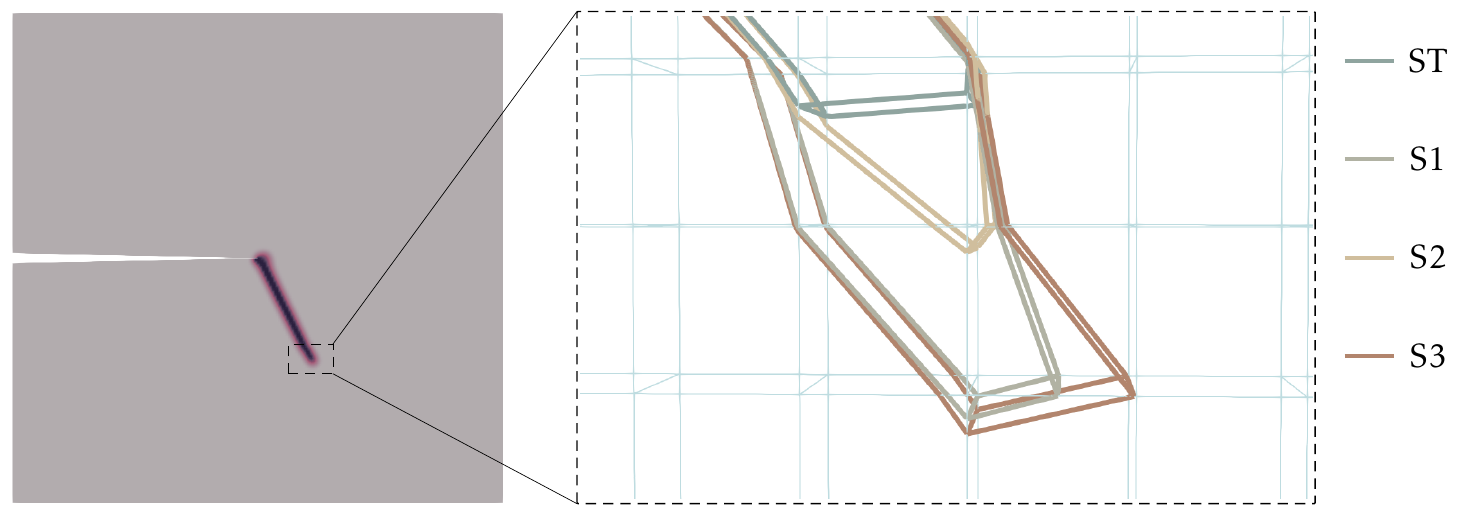}
   \caption{Crack geometry at the final stage in the shear test}
\label{fig_local_global_shear}
\end{figure}

\begin{figure}[!htp]    \centering
   \includegraphics[width=\textwidth]{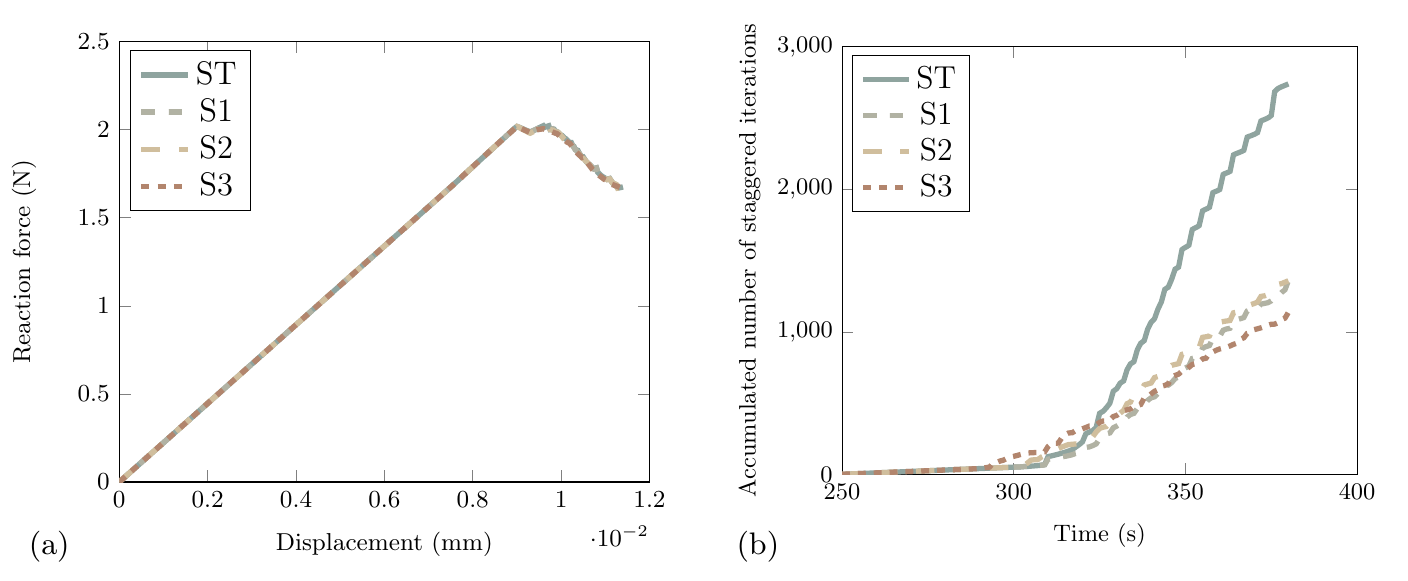}
   \caption{The results of the shear tests: (a) reaction force versus displacement; (b) accumulated number of staggered iterations versus loading time.}
   \label{fig_accu_iteration_shear}
\end{figure}

Fig. \ref{fig_accu_iteration_shear}(a) shows the horizontal reaction force with respect to the displacement, where consistent results are obtained by using different staggered schemes. However, the computational costs of the four staggered schemes are quite different, discussed as follows.
\begin{itemize}
\item[$\blacktriangleright$] The total CPU time corresponding to each of the schemes is ranked as  $T_\text{cpu}^{ST} > T_\text{cpu}^{S1}\approx T_\text{cpu}^{S2} \approx T_\text{cpu}^{S3}$, with the superscript indicating the scheme. The shortest CPU time, i.e., $3.57 \times 10^5 $ s achieved by the S2 scheme, is 66\% of the longest CPU time which is $5.39 \times 10^5 $ s taken when using the ST scheme.
\item[$\blacktriangleright$] The maximum number of iterations among different staggered schemes occurs at the softening stage at different time steps. The maximum number is proportional to the total CPU time, ranked as $N_\text{ST, max}^{ST} > N_\text{ST, max}^{S1}>N_\text{ST, max}^{S2}>N_\text{ST, max}^{S3}$. The maximum numbers of iterations in using the S1, S2 and S3 staggered schemes are 39\%,  42\% and 29\% of $N_\text{ST, max}^{ST}$.
\item[$\blacktriangleright$] Fig. \ref{fig_accu_iteration_shear} compares the total numbers of staggered iterations. It reveals that the S3 scheme requires the fewest staggered iterations, 1138, which is 42\% of the maximum number of staggered iterations, 2739.
\item[$\blacktriangleright$] With the same final displacement, the proposed scheme yields the longest shear crack, cf. Fig. \ref{fig_local_global_shear}. In other words, the crack is developed more rapidly when the proposed scheme is applied, owing to the consideration of the coupling term in solving the evolution equation.
\end{itemize}

\subsection{L-shape panel tests}\label{sec:L-shape}
\begin{figure}[!htp]
   \centering
   \includegraphics[width=0.6\textwidth]{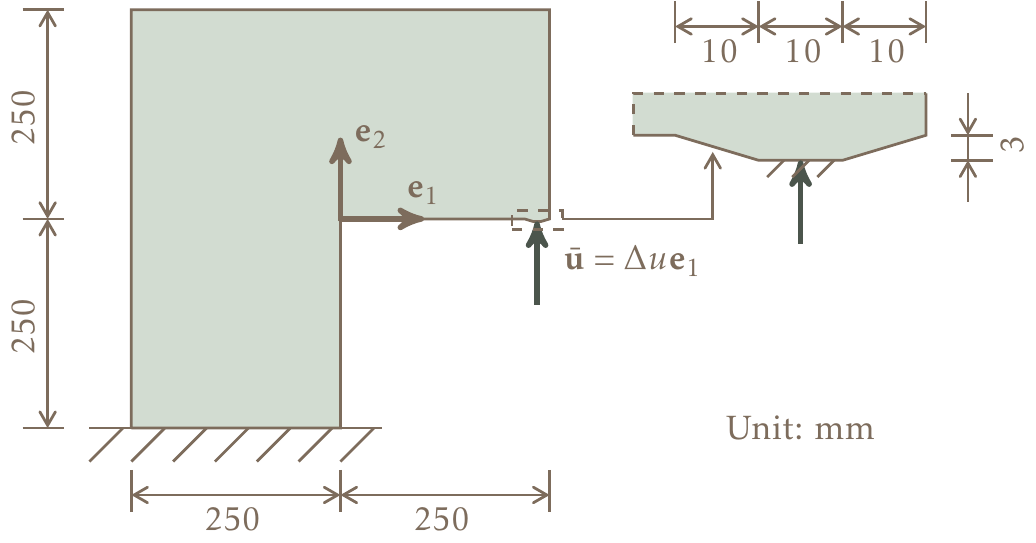}
   \caption{The geometry of the L-shape panel tests and the loading scheme.} 
   \label{geo_bc_L_shape}
\end{figure}

In the last example, we consider an L-shape panel subject to a lifting force at the end of the cantilever. Fig. \ref{geo_bc_L_shape} shows the geometry of the L-shape panel and the loading scheme. This example has been simulated in \cite{gerasimov2016line, wu2020phase,brun2020iterative} and the material parameter setting thereof is used here as given in Table \ref{tab_parameters_L_shape}. 
\begin{table}[htp]
 \centering\small
\vspace{5mm}
\caption{Parameter setting of the L-shape panel}
\begin{tabular}{ccccccc}
  \toprule
&&&&&&\\[-4pt]
Parameter & $\mu$ & $l$ & $\lambda$  & $\eta$ & $G_c$\\[4pt]
\midrule
&&&&&&\\[-4pt]
Value & $1.095 \times 10^{4} \, \text{MPa}$ &  $10\,\text{mm}$&$6.161 \times 10^{3}\, \text{MPa}$ &  $1\times 10^{-6}$ & $0.095 \,\text{N/mm}$& \\[4pt]
\bottomrule
\end{tabular}
 \label{tab_parameters_L_shape}
\end{table}
The right end of the cantilever is subjected to a horizontal displacement at the bottom. The displacement is applied at the bottom of a trapezoid block to avoid stress concentration, as shown in Fig. \ref{geo_bc_L_shape}. Moreover, the phase variable around that area is forced to zero to prevent cracking there. The loading steps are given by,
\be
\Delta u = 1\times 10^{-3} \text{ mm/s}, \qquad \Delta t_n=\left\{
\ba{ll}
0.01, \quad&t \leq 0.2 \text{ s},\\
0.001,\quad& t \leq 0.6 \text{ s}.
\ea
\right.
\ee

\begin{figure}[!htp]
   \centering
   \includegraphics[width=\textwidth]{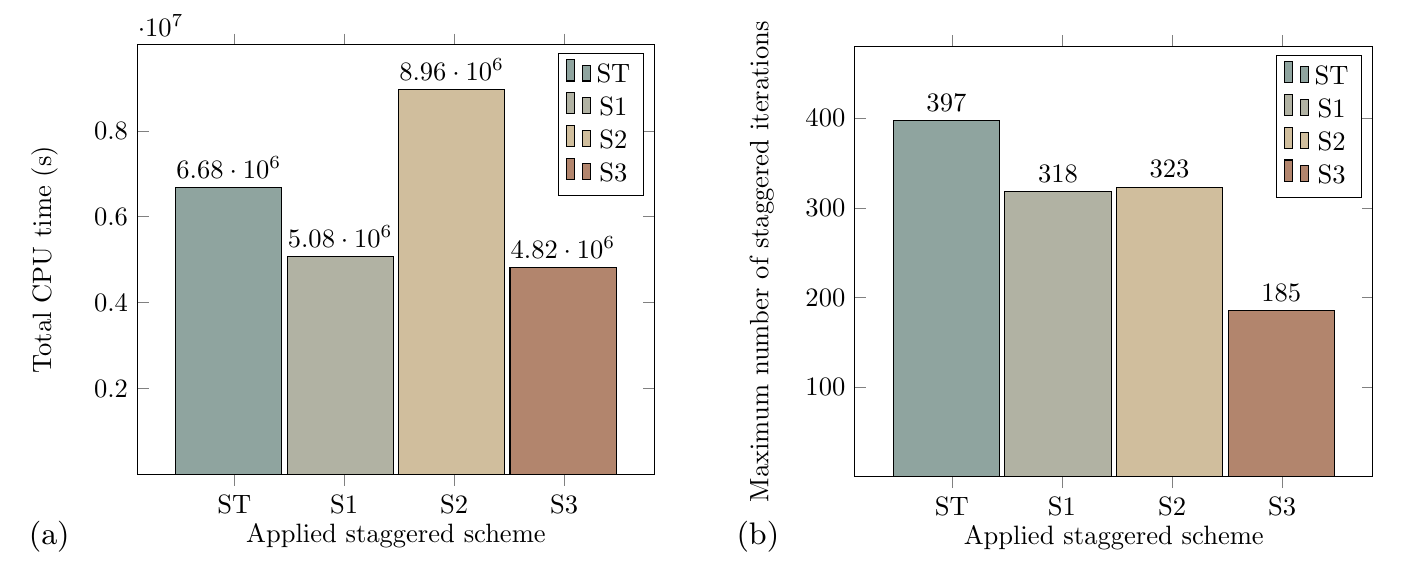}
   \caption{Comparison of the computational costs in the L-shape panel tests: (a) total CPU time; (b) the maximum number of iterations in one timestep.}
 \label{fig_L_shape_summary}
\end{figure}

\begin{figure}[!htp]
    \centering
    \includegraphics[width=\textwidth]{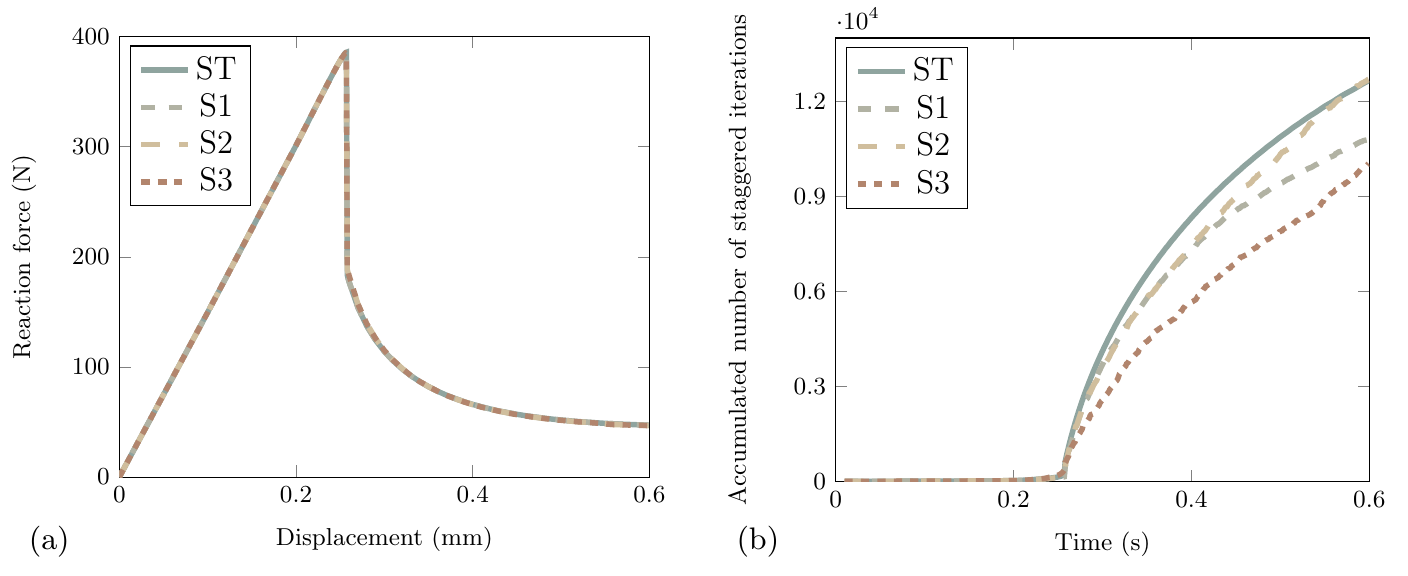}
    \caption{The results of the L-shape panel tests: (a) reaction force versus displacement; (b) accumulated number of staggered iterations versus loading steps.}
    \label{fig_disp_accu_iter_L_shape}
 \end{figure}

\begin{figure}[!htp] 
   \centering
   \includegraphics[width=0.8\textwidth]{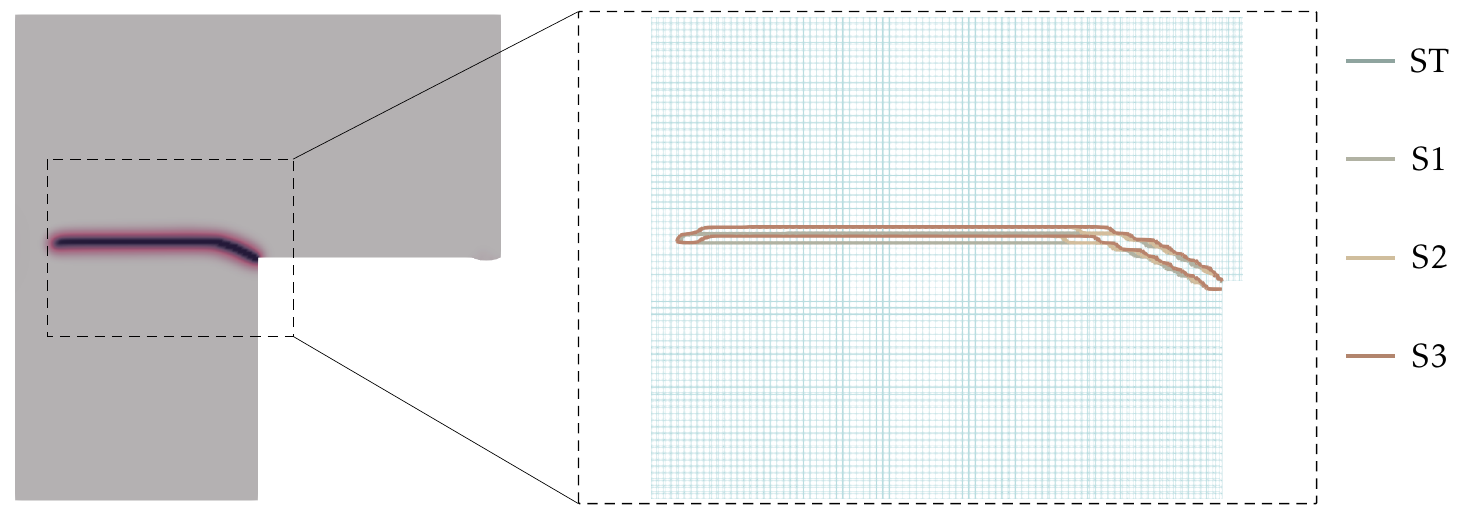}
   \caption{A comparison of the crack geometries at the final stage in the L-shape panel tests computed using different staggered schemes.}
\label{fig_local_global_L_shape}
\end{figure}

\begin{figure}[!htp]
   \centering
   \includegraphics[width=\textwidth]{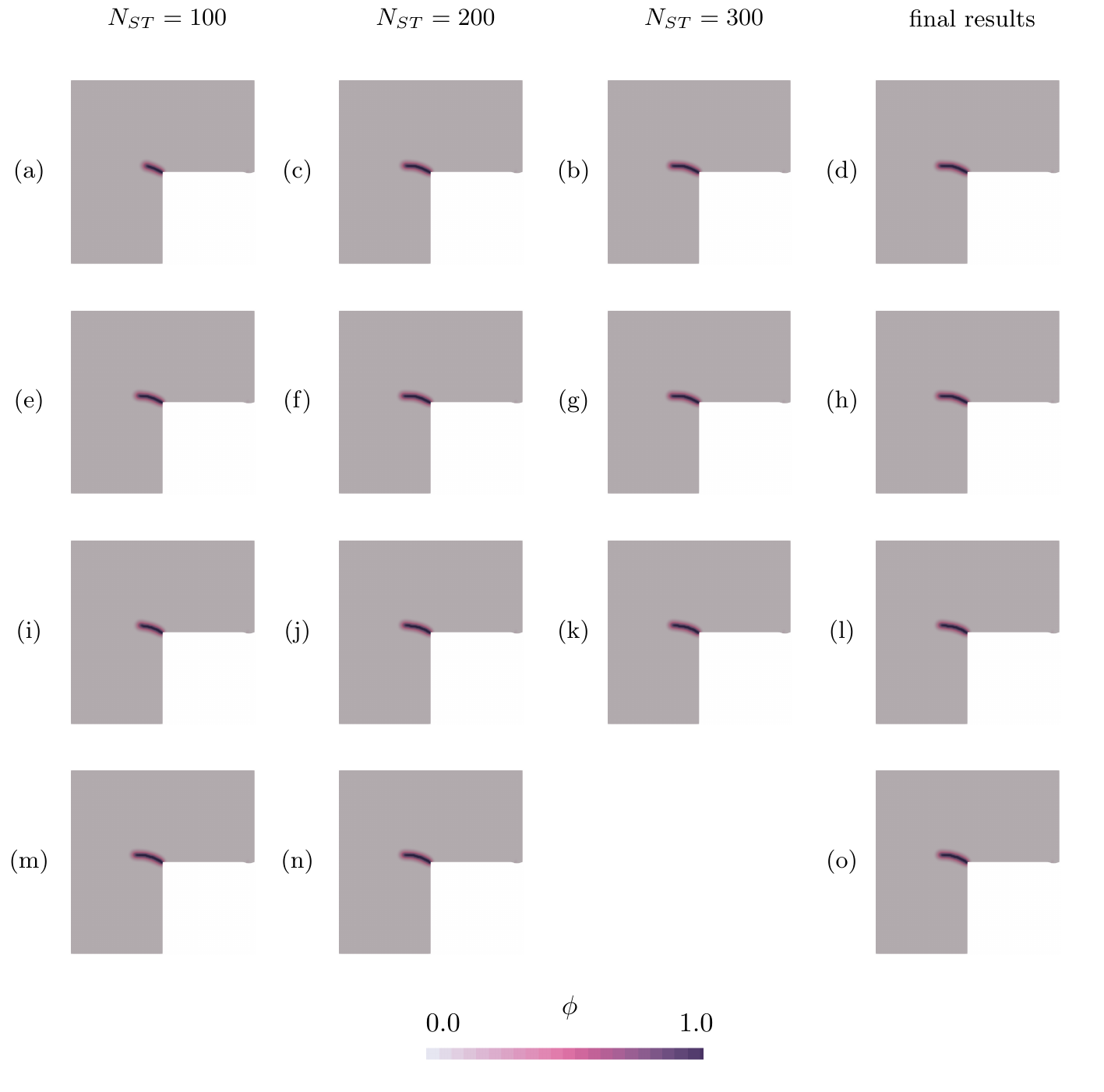}
   \caption{Crack propagation within the staggered iterations at time 0.258 s in the L-shape panel test when using different schemes: (a-d) ST scheme;(e-h) S1 scheme; (i-l) S2 scheme; (m-o) S3 scheme. }
\label{fig_crack_propagation_L_shape}
\end{figure}

Fig.~\ref{fig_disp_accu_iter_L_shape}(a) exhibits the relation between the displacement and reaction force on the tip of the cantilever, which shows the consistency of the results derived by four staggered schemes. Referring to the crack patterns provided in Fig.~\ref{fig_local_global_L_shape}, the S2 and S3 schemes yield a crack slightly bended upwards, and the lengths are very close. Although simulation results are similar, four schemes achieve different computational efficiency:
\begin{itemize}
\item[$\blacktriangleright$] The relation of the total CPU time is distinct from the previous examples, which is  $T_\text{cpu}^{S2}  > T_\text{cpu}^{ST} > T_\text{cpu}^{S1}> T_\text{cpu}^{S3}$. The CPU time using the S2 scheme is 1.76 times of that corresponding to the ST scheme. Moreover, the S3 scheme consumes the least time, $4.82 \times 10^6$ s, which is 72 \% of the one using the ST scheme.
\item[$\blacktriangleright$]  In terms of the maximum number of staggered iterations, the present example yields a similar conclusion as the previous two examples, e.g., the maximum number of staggered iterations when using the proposed scheme is generally smaller than the one using the ST scheme. For the current example, the maximum numbers of iterations respectively using the S1, S2 and S3 staggered schemes are 80\%,  81\% and 46\% of $N_\text{ST, max}^{ST}$.
\item[$\blacktriangleright$] The total staggered iterations are displayed in Fig.~\ref{fig_disp_accu_iter_L_shape}(b). It uncovers that the S3 scheme has the fewest staggered iterations, 2988, which is 79\% of the ST staggered iterations, 12645. Comparing the curves in the Fig. \ref{fig_disp_accu_iter_L_shape}(b), we find that most reduction occurs at the crack nucleation and the early stage of crack propagation. Fig. \ref{fig_crack_propagation_L_shape} demonstrates that a crack of the same length is formed within fewer iterations when using the proposed schemes.
\end{itemize}
The reason scheme S2 fails in reducing the computation time can be explained by its assumption as given in Table \ref{tab_summary_staggered_scheme}. The S2 scheme guarantees that the incremental of the second invariant is zero, provided that the strain trace is zero and the resultant evolution of the first invariant is automatically zero and thus neglected. However, in the present example, this assumption is violated, as dilation occurs around the crack. Moreover, unlike the previous two examples, the violation of the assumption cannot be corrected by solving the momentum balance as the top boundary is traction-free.

\subsection{Summary of numerical tests}
Based on the above results, we conclude that
\begin{itemize}
\item By using the proposed schemes, the derived results are consistent with the one obtained using the standard scheme;
\item The proposed schemes can accelerate the convergence, e.g., reduce the maximum number of staggered iterations;
\item With a higher convergence rate, generally speaking, the proposed schemes are computationally more efficient as they cost less CPU time and require fewer staggered iterations in total; 
\item The schemes S1 and S3 which fix first and second invariants of the stress tensor unconditionally accelerate the numerical computation, while the scheme S2 improves computational efficiency under the condition of a negative strain trace. When this condition is violated, the S2 scheme can be less efficient than the ST scheme.
\end{itemize}
Apart from the staggered iteration, we also notice that both the momentum balance and the evolution equation are nonlinear and solved in an iterative manner. In Table \ref{tab_summary_wall_clock}, the accumulated number of NR iterations for solving displacement and the phase variable, together with the wall-clock time, are compared. With regard to computational performance of the ST and S3 schemes for the shear tests, we discuss the tradeoff as follows
\begin{itemize}
\item The NR iteration for solving the evolution equation is much cheaper than that for the momentum balance in terms of computational cost since $d$ is a scalar while $\mathbf{u}$ is a vector. In the shear test,  the average wall-clock time or cpu time for solving the evolution equation, 0.047 s, is roughly 10\% of that for solving the momentum balance, 0.490 s.
\item As the proposed staggered schemes update the coupling term, i.e., the active strain energy, when solving the evolution equation, the total number of NR iterations for solving $d$, 8141, is usually larger than that in using the ST scheme, 4882.
\item Thanks to the reduced total number of staggered iterations in the proposed scheme, the total number of NR iterations for solving $\mathbf{u}$, 2563, is considerably reduced as compared to the ST scheme, 4021.
\item A comparison of the percentages of the wall-clock time clearly shows that the fast schemes spend more time, 360 s, in solving the evolution equation than the ST scheme, 228 s. However, solving the evolution equation takes much less time as compared to solving the momentum balance equation. In other words, the increased wall clock time in solving $d$, 132 s, is incomparable to that in solving $\mathbf{u}$, 780 s. Therefore, the total wall clock time is reduced by applying the proposed schemes.
\end{itemize}

\begin{table}[htbp]
    \small
    \caption{Summary of computational costs based on the standard and serendipity staggered schemes.}
    \label{tab_summary_wall_clock}
     \centering\small
     \begin{tabular*}{\textwidth}
       {@{\extracolsep\fill}ccccccccccccc@{\extracolsep\fill}}
    \toprule
    \multirow{2}{*}{Test cases}  && \multicolumn{2}{c}{ST} && \multicolumn{2}{c}{S1}
    && \multicolumn{2}{c}{S2} && \multicolumn{2}{c}{S3} \\
    \cline{3-4} \cline{6-7} \cline{9-10} \cline{12-13}\\[4pt]
    &&$d$&$\bf u$& &$d$&$\bf u$& &$d$&$\bf u$& &$d$&$\bf u$\\
    \midrule
    \multirow{5}{*}{Tensile test}
    &$N_{NR}$                      &1654&881    &&1418&656    &&1826&599    &&1648&602\\[5pt] 
    &$T_{NR}$ (s)                  &79.6&415    &&74.6&344    &&105&304     &&35.5&264 \\[5pt] 
    &$\bar{T}_{NR}$ (s)            &0.048&0.471 &&0.053&0.532 &&0.058&0.508 &&0.022&0.438\\[5pt] 
    &$T_{NR} / T_{NR}^\text{total}$&16\%&84\%   &&18\%&82\%   &&24\%&76\%   &&26\%&74\%\\[5pt]  
    \midrule
    \multirow{5}{*}{Shear test}
    &$N_{NR}$                      &4882&4021   &&7559&2568   &&7180&2452   &&8141&2563\\[5pt] 
    &$T_{NR}$ (s)                  &228&1970    &&334&1160    &&332&1160    &&360&1190\\[5pt] 
    &$\bar{T}_{NR}$ (s)            &0.047&0.490 &&0.044&0.452 &&0.046&0.473 &&0.044&0.464\\[5pt] 
    &$T_{NR} / T_{NR}^\text{total}$&10\%&90\%   &&22\%&78\%   &&22\%&78\%   &&23\%&77\%\\[5pt]  
   \midrule
    \multirow{5}{*}{L-shape panel}
    &$N_{NR}$                      &20681&18500 &&49941&11600 &&87408&19700 &&38470&11600\\[5pt]
    &$T_{NR}$ (s)                  &2850&24477  &&3700&21170  &&9650&25104  &&3150&19904\\[5pt]
    &$\bar{T}_{NR}$ (s)            &0.138&0.756 &&0.074&0.548 &&0.110&0.785 &&0.081&0.583\\[5pt]
    &$T_{NR} / T_{NR}^\text{total}$&13 \%&87\%  &&24\%&76\%   &&33\%&67\%   &&21\%&79\%\\[5pt]
  \bottomrule 
 \end{tabular*}
\end{table}

\section{Conclusion}
We propose novel staggered schemes to reduce the number of staggered iterations and improve the convergence rate in phase-field modeling of brittle fractures. Inspired by the fixed-stress staggered scheme for solving coupled problems in poromechanics, we propose to fix principle invariants of the stress in solving the evolution equation instead of fixing displacement in the standard staggered scheme for phase-field fracture model. Then we are able to derive the relation between the increments of the strain and the phase variable. By using this relation, the phase variable can be solved from the evolution equation, with updated coupling term in the equation, namely the active strain energy. Subsequently, the momentum balance equation is solved. As compared to the standard scheme, the critical difference lies in that we update the active strain energy based on its first-order coupling relation, namely $\Delta \bep$ and $\Delta d$. More precisely, we predict the change of the active strain energy induced by $\Delta d$ and solve the evolution equation based on this prediction. Since the active strain energy not only depends on the displacement derived from the momentum balance but also depends on $\Delta d$, the required number of staggered iterations is decreased. On the contrary, in the standard staggered scheme, the $\Delta d$-induced change of the displacement and the active strain energy can only be updated after solving the momentum balance.

Under the proposed concept, we further present three staggered schemes by fixing the first or the second or both invariants of the stress. The schemes are verified by comparing with the ST scheme in three typical benchmark tests. The results show that the S1 and S2 schemes can achieve a dramatic reduction of the maximum number of staggered iterations and the total CPU time in all three numerical experiments. The S2 scheme performs comparably except that in the shear test where the underlying assumption of a negative strain trace is violated and hence the computational efficiency is not improved. A comparison of the crack propagation during the staggered iteration within one time step reveals that the crack obtained by the proposed staggered scheme propagates more rapidly owing to the update of the active strain energy. 

Beside the achieved improvement on computational efficiency, our schemes are also featured by a clear physical interpretation and no modification of the original equations. Moreover, no hyper parameters such as the over-relaxation parameter or stabilization parameters are introduced in the numerical computation. The fast staggered schemes can be conveniently implemented by adding two modules, respectively for computing the stiffness component $\tilde{\bf K}_{d\bep}$ and updating the active strain energy, to the standard staggered scheme. As the S1 and S3 schemes are valid without any prerequisites, they are recommended for applications.

\appendix
\section{Supplementary deduction} \label{summplementry_deduction}
In this appendix, the deduction of Eq. \eqref{eq:trdep}, Eq. \eqref{eq:epdep} and Eq. \eqref{e_dev_1/2} is provided. To expand Eq. \eqref{eq_assumption1}, we first recall the following rules of tensor algebra
\begin{equation}
\left[
    (\bI\otimes\bI)^{\overset{24}{T}}\mathbf{A}
\right]\cdot \bI =\tr \mathbf{A}, \quad \left[(\bI\otimes\bI)\mathbf{A}
\right]\cdot \bI =3 \, \tr \mathbf{A}.
\end{equation}
Insertion of the constitutive equation \eqref{eq_sigma} into Eq. \eqref{eq_assumption1} yield the following relation
\begin{equation}
    g(d)\left[
        2\mu(\tr \Delta \bep)-\frac{2\mu}{3}(3\tr \Delta \bep)
        +k\bI^+(3\tr \Delta \bep)
    \right]
        +k\bI^-(3\tr \Delta \bep)
    -2(1-d)\langle3\tr \Delta \bep\rangle_+]\Delta d=0.
\end{equation}
The preceding equation can be rewritten in a compact form
\begin{equation}
    \begin{split}
        [g(d)\bI^++\bI^-]\tr \Delta \bep
        &=2(1-d)\langle\tr \bep\rangle_+\Delta d.\\
    \end{split}
\end{equation}
After taking the definition from Eq. \eqref{eq_definition_+}, Eq. \eqref{eq:trdep} is derived.

Before turning to \eqref{eq_assumption2}, we recall the following relations
\begin{equation}
    \left[
        (\bI\otimes\bI)^{\overset{24}{T}}
    \right]\cdot\mathbf{B}
    =\mathbf{A}\cdot\mathbf{B}, \quad
        \left[
       ( \bI\otimes\bI) \mathbf{A}
    \right]\cdot \mathbf{B}
    =(\tr \mathbf{A})(\tr \mathbf B).
\end{equation}
Using these relations along with $\mathbf{A}^\text{dev} \cdot \bI = 0$, Eq.  \eqref{eq_assumption2} is rewritten as
\begin{equation}
    g(d)2\mu\bep^{\dev}\cdot
    \Delta \bep-
    2(1-d)(2\mu\bep^{\dev}\cdot \bep^{\dev})\Delta d=0,
\end{equation}
which directly results in Eq. \eqref{eq:epdep}. To derive Eq. \eqref{e_dev_1/2}, we firstly expand the left-hand side of the equation as
\begin{equation}\label{eq:epn1epn1}
    \bep^{\dev,k+\frac{1}{2}}\cdot\bep^{\dev,k+\frac{1}{2}}=(\bep^{\dev,k}+\Delta \bep^{\dev})\cdot(\bep^{\dev,k}+\Delta \bep^{\dev})
    =\bep^{\dev,k}\cdot \bep^{\dev,k}
    +2\bep^{\dev,k}\cdot \Delta \bep^{\dev}
    +\Delta \bep^{\dev}\cdot
    \Delta \bep^{\dev}.
\end{equation}
The first term is already known and the second term can be computed from Eq. \eqref{eq:epdep}. Then the only unknown term $ \Delta \bep^{\dev}\cdot
    \Delta \bep^{\dev}$ can be rewritten by inserting the definition of a deviatoric tensor, $\bep^{\dev} = \bep -\frac{1}{3}(\tr \bep^k)\bI$, as
\begin{equation} \label{eq_dot_de_dev}
    \Delta \bep^{\dev}\cdot
    \Delta \bep^{\dev}
    =\left[\Delta \bep -\frac{1}{3}(\tr \Delta \bep)\bI\right]
    \left[\Delta \bep -\frac{1}{3}(\tr \Delta \bep)\bI\right]
    =\Delta \bep\cdot\Delta\bep-
    \frac{1}{3}(\tr\Delta \bep)^2.
\end{equation}
From Eq.~\eqref{eq:trdep},  this dot production equals zero when the strain trace is negative,
\begin{equation}
    \tr \bep < 0 \rightarrow \tr \Delta\bep =  0 \rightarrow  \Delta \bep^{\dev}\cdot
    \Delta \bep^{\dev}= 0.
\end{equation}
Then we discuss the results when the strain trace is positive. We assume the dot production between the incremental stress and strain as
\begin{equation}
    \Delta \bsg \cdot \Delta \bep = 0,
\end{equation}
since $\Delta \bsg = \mathbf{0}$. Inserting Eq. \eqref{eq_delta_stress} into the above equation results in
\begin{equation}\label{eq:gd}
    \begin{split}
    g(d)\left[
        2\mu\Delta\bep\cdot\Delta\bep-\frac{2}{3}\mu(\tr \Delta\bep)^2+k\bI_+
        (\tr\Delta\bep)^2
    \right]+k\bI_-(\tr\Delta\bep)^2\\
    -2(1-d)\left[
        2\mu\bep^{\dev}\cdot\Delta\bep
        +k\langle\tr\bep\rangle_+(\tr \Delta\bep)
    \right]\Delta d=0.
    \end{split}
\end{equation}
Noting that the strain trace is positive, the following relation holds
\begin{equation}
    \tr \bep>0, \, \tr \Delta\bep\ne 0, \, \bI_- = \mathbf{0},
\end{equation}
based on which, Eq. \eqref{eq:gd} is rewritten in a compact form
\begin{equation}
    g(d)[2\mu\Delta\bep\cdot\Delta\bep+
    \lambda(\tr\bep)^2]-2(1-d)
    [2\mu\bep^{\dev}\Delta\bep+
    k\langle\tr\Delta\bep\rangle_+(\tr\Delta\bep)]\Delta d=0.
\end{equation}
From the preceding relation, we can derive the dot product of the incremental strains and themselves as
\begin{equation} \label{dot_increment_strain_strain}
    \Delta\bep\cdot\Delta\bep=-\frac{\lambda}{2\mu}(\tr \Delta\bep)^2+\frac{2(1-d)}{g(d)}
    \left[
        \bep^{\dev}\cdot\Delta\bep
        +\frac{k}{2\mu}\langle\tr \bep\rangle_+
        \tr\Delta \bep
    \right]\Delta d.
\end{equation}
Insertion of Eq. \eqref{dot_increment_strain_strain}, Eq. \eqref{eq:trdep} and Eq. \eqref{eq:epdep} into Eq. \eqref{eq_dot_de_dev} yields
\begin{equation} \label{eq_dot_de_dev1}
    \Delta\bep^\text{dev} \cdot\Delta\bep^\text{dev}= 4\frac{(1-d)^2}{g^2(d)}
    (\bep^{\dev} \cdot \bep^{\dev}) \Delta d.
\end{equation}
Finally, with Eq. \eqref{eq:epdep} and Eq. \eqref{eq_dot_de_dev1} at hand, we reach Eq. \eqref{e_dev_1/2}.

\section{Reinterpretation of the fixed-stress splitting scheme in poromechanics} \label{reinterpretation_fixed_stress}
In this section, we follow the idea of Eq. \eqref{eq_no_idea_for_name} to explain the underlying mechanism of the fixed-stress splitting scheme from an alternative point of view. First, we recall the governing equations for a typical poromechanical problem, cf. \cite{kim2011stability1},
\be
\left\{
\ba{lcl}
\dot{m} + \div \mathbf{w} &=& \rho_{f,0} f,\\
\div \bsg + \rho \mathbf{g} &=& \mathbf{0},\\
\ea
\right.
\ee
where $m$ is the fluid mass per unit bulk volume, $\mathbf{w} $ the fluid mass flux, $ \rho_f$ the initial fluid density. The total stress is composed of the effective stress, $\bsg'$, and the fluid pressure, $p$, as
\be
\bsg = \bsg'- b p \mathbf{I},
\ee
where $b$ is the Biot coefficient. Provided that the solid is a linear elastic material, the effective stress can be computed from the strain as
\be
\bsg'= 2\mu \bep^\text{dev} + k (\tr \bep) \mathbf{I}.
\ee
The total density $\rho$ is defined as $\rho = n^\rho_f + (1-n)\rho_s$, where $\rho_s$ is the solid density, $n$ the so-called true porosity. In order to solve this system in a sequential manner, we need to solve the first equation under the assumption that the coupling term in the second equation is unchanged. Therein, the coupling term in the momentum balance is the total stress. Following Eq. \eqref{eq_no_idea_for_name}, we first assume the stress trace is constant and yield the following relation
\be \label{fixed_stress_relation}
\frac{\text{d}}{\text{d}\, t}I_1^{\bsg} = 3 k \tr  \dot{\bep} - 3 b \dot{p} = 0 \rightarrow \tr \dot{ \bep} = \frac{b}{k} \dot{p}.
\ee
Then we recall the constitutive equation for $m$
\be
\dot{m} = b \tr \dot{\bep} + \frac{1}{M}\dot{p},
\ee
where $M$ is the Biot modulus. By inserting Eq. \eqref{fixed_stress_relation} into this constitutive equation, we derive
\be
\dot{m} =\left( \frac{b^2}{k} + \frac{1}{M}\right) \dot{ p}.
\ee
Substituting $\dot{m}$ in the first equation with the above relation, we finally obtain the form which is used for the fixed-stress splitting scheme. Now we can further let the second invariant also be a constant. Then we derive
\be \label{fixed_stress_relataion1}
\frac{\text{d}}{\text{d}\, t}I_2^{\bsg} = \frac{\text{d}}{\text{d}\, t}I_2^{\bsg'} = 0.
\ee
As the pressure term is an isotropic tensor, its second invariant vanishes. Therefore, assuming the second invariant of the total stress fails in deriving any further relation between the pressure the strain. 

\section{A one-dimensional example}
\begin{figure}[!htp]
   \centering
   \includegraphics[width=\textwidth]{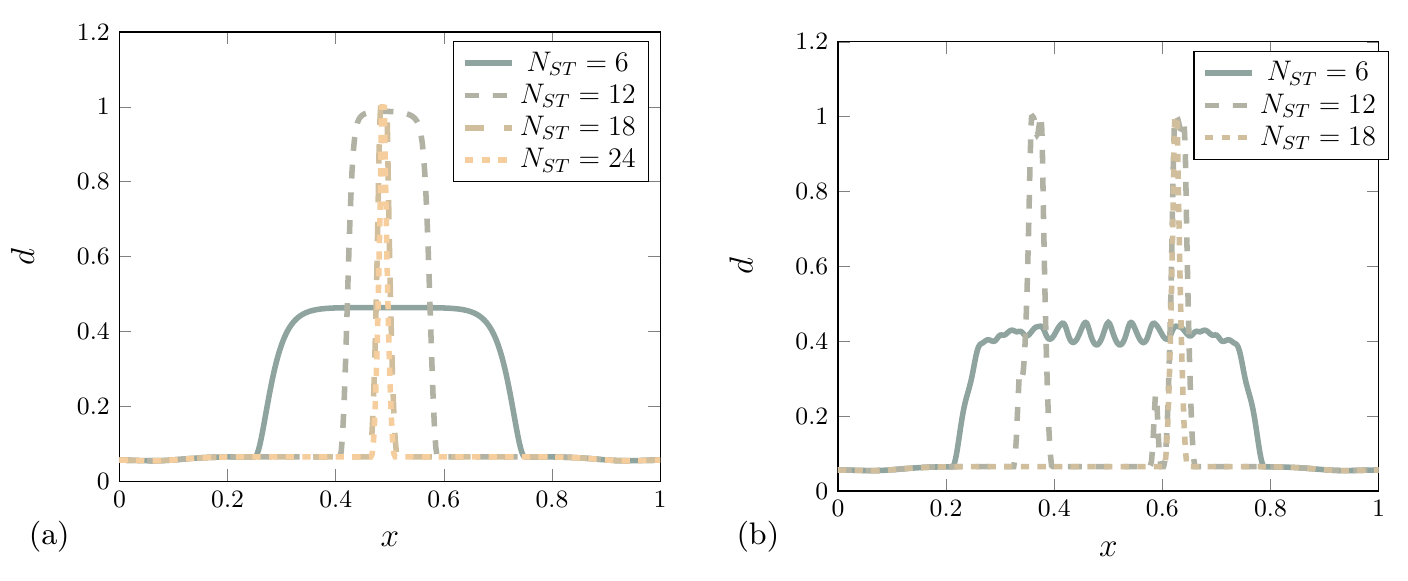}
   \caption{The evolution of the phase variable in the one-dimensional example: (a) using the ST scheme; (b) using the proposed staggered scheme.}
   \label{phi_1d}
\end{figure}

Here, a simple one-dimensional example is provided to demonstrate how the proposed scheme accelerate the crack nucleation. We first rewrite the evolution equation \eqref{eq_evolution_equation} as
\be
-(1-d) E \varepsilon^2 + \frac{G_c}{c_w}\left(\frac{1}{l} + 2l d''\right) = 0,
\ee
where $(\ )^\prime$ is the spatial derivative operator. Following the assumption that the stress is constant when solving the evolution equation, we derive
\be
\frac{\text{d}}{\text{d}t}\sigma = \frac{\text{d}}{\text{d}t}(1-d)^2 E\varepsilon = 0 \rightarrow \Delta \varepsilon = \frac{2 \varepsilon}{1-d} \Delta d.
\ee
After solving the evolution equation, we update the strain trace by $\Delta d$ as
\be
\varepsilon^{k+\frac{1}{2}} = \varepsilon^k + \frac{2 \varepsilon}{1-d} \Delta d.
\ee
Correspondingly, the active strain energy will be updated as
\be
\Psi^{k+\frac{1}{2}} = \Psi ^k + \frac{4E\varepsilon}{1-d} \Delta d + \frac{4E\varepsilon^2}{(1-d)^2} \Delta d ^2.
\ee
If the condition in Eq.~\eqref{condition_update} is satisfied, $\Psi^{k+\frac{1}{2}}$ must be greater than $\Psi ^k $, which then leads to further evolution of the phase variable and accelerate the crack nucleation. Herein, one dimensional example with a length of 1 is set up, where $E=1$, $G_c = 0.1 $ and $l = 0.01$. The sample is extended at one end by $\Delta \bar{u} = 0.01 $ and fixed at the other end. The evolution of the phase variable is accomplished within one time step, at time $t=201$, where the standard staggered scheme needs 25 staggered iterations while the proposed scheme needs 17. The evolution along the example with respect to different staggered iterations are given in Fig. \ref{phi_1d}. We observe that in the early stage the distribution of the phase variable in the possible cracking region oscillates owing to the local update of the active strain energy. This fluctuation helps to accelerate the localization. 

\section*{Declaration of competing interest}
The authors declare that they have no known competing financial interests or personal relationships that could have appeared to influence the work reported in this paper.

\bibliographystyle{elsarticle-num}
\bibliography{fast_staggered_scheme_luo.bib}
\end{document}